\documentclass[11pt]{article}
\usepackage{amsmath,comment}
\usepackage{mathrsfs}
\usepackage{amsfonts}

\usepackage{amsfonts, amsmath, amssymb}
\usepackage{amssymb,amsfonts,amsmath,color,
latexsym, epsfig,cite, psfrag,eepic,colordvi}
\usepackage{amscd,graphics}
\usepackage{lineno}

\newcommand{\qed}{\hfill $\Box $}
\newcommand{\pf}{\noindent {\bf Proof.} }

\newtheorem{theorem}{Theorem}[section]
\newtheorem{lemma}[theorem]{Lemma}

\newtheorem{myDef}[theorem]{Definition}
\newcommand{\ex}{{\rm ex}}
\newcommand{\G}{{\cal G}}

\begin{document}

\title{Rainbow Matchings in $K_{r+1}$-free Graphs}
\author{
Hongliang Lu$^{1}$\thanks{
This work was supported by National Key Research and Development Program of China  2023YFA1010203
and  the National Natural Science Foundation of China under grant No.12271425
},
Xinxin Ma$^{1}$\thanks{
Corresponding author. E-mail:
\texttt{maxinxinzzu@163.com}.
},
and Xing Peng$^{2}$
\\[0.5em]
{\small $^{1}$School of Mathematics and Statistics,
Xi'an Jiaotong University, Xi'an, Shaanxi 710049, China}
\\
{\small $^{2}$School of Mathematics and Statistics,
Anhui University, Hefei, Anhui 230601, China}
}

\date{}

\maketitle

\begin{abstract}
 Let $\G=\{G_1,\ldots,G_s\}$ be a system of graphs on the common vertex set with $n$ vertices. A graph $F$ with $s$ edges is rainbow if there is a bijection $\phi: E(F) \to [s]$ such that $e \in E(G_{\phi(e)})$ for each $e \in E(F)$.
For a family of graphs $\cal F$, the Tur\'an number $\ex(n, \cal F)$ is the maximum number of edges in an $n$-vertex graph which does not contain a member of $\cal F$ as a subgraph.
Akiyama and Frankl proved that for $n \geq 2s$, the condition  $e(G_i) > \ex(n,M_s)$ for each $i \in [s]$ guarantees  a rainbow matching $M_s$. In this paper, we establish an analogous result for graph systems with bounded clique number.
Namely, for positive integers $n,r,m$ with $n\geq 2m+2$ and $n$ sufficiently large in terms of $r$, we prove that if  $G_i$ is $K_{r+1}$-free and $e(G_i)> \ex(n,\{K_{r+1},M_{m+1}\})$ for each $i \in [m+1]$,
then $\{G_1,\ldots,G_{m+1}\}$ admits a rainbow matching $M_{m+1}$.
This result can be viewed as a rainbow analogue of a theorem of Alon and Frankl on the maximum number of edges in graphs with bounded clique number and matching number  for sufficiently large $n$.


\end{abstract}

\section{Introduction}
A typical problem in extremal graph theory is to determine  optimal sufficient conditions for a graph  to contain  a prescribed graph  as a subgraph. The condition in terms of the number of edges is an extensively-studied one.
 For a family of graphs $\cal F$, the Tur\'an number $\ex(n, \cal F)$ is the maximum number of edges in an $\cal F$-free $n$-vertex graph. Here a graph is $\cal F$-free if it does not contain any member of $\cal F$ as a subgraph.
 For the case where $\cal F$ consists of a single graph $F$, we write $\ex(n,F)$ instead.
  Thus every $n$-vertex graph $G$ with $e(G)>\ex(n, F)$ contains $F$ as a subgraph. 
 The first result on the Tur\'an number of a graph perhaps is Mantel's result on $\ex(n,K_3)$. The well-known Tur\'an's theorem gives the Tur\'an number of $K_{r+1}$ for each $r \geq 2$. A classical result concerns matching. A matching $M_{m+1}$ is a set of $m+1$ pairwise vertex-disjoint edges. Erd\H{o}s and Gallai proved that for all $m \geq 1$ and $n \geq 2m+2$,
\[
\ex(n,M_{m+1})=\max\left\{\binom{2m+1}{2},\binom{m}{2}+m(n-m-1) \right\}.
\]
.

A natural question is  whether sharp sufficient conditions for subgraph containment extend to the rainbow setting.  Let $\G=\{G_1,\ldots,G_s\}$ be a system of graphs on the common $n$-vertex set $V$. A graph $F$ with $s$ edges is {\it rainbow} if there is a bijection $\phi: E(F) \to [s]$ such that $e \in E(G_{\phi(e)})$ for each $e \in E(F)$. Here $[i]=\{1,2,\ldots,i\}$ for a positive integer $i$. In particular, one may ask whether an edge
condition that forces a copy of $F$ in a single graph also forces a rainbow
copy of $F$ when imposed on every graph in the system.
 A pioneering result in this direction is due to Akiyama and Frankl and can be stated as follows.
\begin{theorem} [Akiyama and Frankl] \label{AkF}
Let $n$ and $m$ be positive integers with $n \geq  2m + 2$ and
let $\G =\{G_1,\ldots, G_{m+1}\}$ be a system of graphs on the common vertex set with $n$ vertices. If
\[
e(G_i)>\ex(n,M_{m+1})
\]
 for each $i \in [m+1],$ then $\G$ admits a rainbow matching $M_{m+1}$.
\end{theorem}
Note that if we take $G_i=G$ for each $i \in [m+1]$, then it recovers the Erd\H{o}s--Gallai theorem.  Thus the same edge condition that guarantees a matching also guarantees a rainbow matching.
There are several other results exhibiting this phenomenon; see, for example \cite{AH17,HLF,LZW,YZZ,Zv}. These results motivate the question of whether analogous
sharp rainbow results hold if we impose additional  restrictions on
graphs in the system. In this paper, we consider graph systems with bounded
clique number.

 Let $T(n,r)$ be the $r$-partite Tur\'an graph which is the balanced complete $r$-partite graph on $n$ vertices. For $r \geq 3$, let $G(n,m,r)$ denote the complete $r$-partite $n$-vertex graph with one part of order $n-m$  and each other part of order $\lceil m/(r-1)\rceil$ or $\lfloor m/(r-1)\rfloor$.
 Alon and Frankl \cite{AF}, and independently  Fu, Wang, and Yang \cite{FWY}, established the following result, which determines the maximum number of edges in a graph with bounded clique number and matching number.
 \begin{theorem} \label{AF}
 Let $n,m,r$ be positive integers  with $r \geq 2$ and $n \geq 2m+2$.
If $G$ is an $n$-vertex $K_{r+1}$-free graph satisfying  $e(G)>\max\{e(G(n,m,r)),e(T(2m+1,r))\}$, then $G$ contains a matching $K_{m+1}$ and this bound is sharp.
 \end{theorem}
Note that  Theorem \ref{AF} is equivalent to $$\ex(n,\{K_{r+1},M_{m+1}\})=\max\{e(G(n,m,r)),e(T(2m+1,r))\}.$$
 For extensions and generalizations of Theorem \ref{AF}, readers can refer to \cite{FWY,GTZ,Ge24,LZL,XK,YYP}.
In this paper, we will prove the following rainbow analogue of Theorem \ref{AF} for large $n$ in terms of $r$.
\begin{theorem}\label{main}
Let $n,r,m$ be three positive integers such that $n\geq 2m+2$ and $1/n\ll 1/r$.  Assume
$\G=\{G_1,\ldots,G_{m+1}\}$ is a system of graphs on the common vertex set with $n$ vertices.
  If  for $1\leq i\leq m+1$, the graph $G_i$ is  $K_{r+1}$-free and
\[
e(G_i)> \max\{e(G(n,m,r)),e(T(2m+1,r))\}=\ex(n,\{K_{r+1},M_{m+1}\}),
\]
then $\G$  admits a rainbow matching $M_{m+1}$.
\end{theorem}
Therefore, for sufficiently large $n$, the optimal edge condition
for a matching in a $K_{r+1}$-free graph remains sufficient in the rainbow
setting.

Since the proof of Theorem 1.3 relies on results concerning hypergraph matchings, we first introduce some related results on hypergraph matchings. Let $k \ge 2$ be a positive integer.
For a set $S$ with at least $k$ elements, let ${S\choose
  k}=\{e\subseteq S\ | \ |e|=k\}$.
A \emph{$k$-uniform hypergraph} or simply a \textit{$k$-graph}  is a pair
$H=(V,E)$, where $V=V(H)$ is a finite set of {\it vertices} and
 $E=E(H)\subseteq {V\choose k}$ is the set of {\it edges}.

A \emph{matching} in a hypergraph $H$ is a subset of $E(H)$ which consists of pairwise
vertex-disjoint edges. The \textit{matching number of} a hypergraph $H$, denoted $\nu(H)$, is  the maximum cardinality of a matching in $H$.
Erd\H{o}s
\cite{Erdos65} in 1965 made the following conjecture:
For positive integers $k,n,t$ with $n \ge kt$, every $k$-graph $H$ on $n$ vertices with $\nu(H) <
t$ satisfies $e(H)\leq \max \left\{ {n\choose k}-{n-t+1\choose k}, {kt-1\choose k} \right\}.$
This bound is tight for the complete $k$-graph on $kt-1$ vertices and for the $k$-graph on $n$ vertices in which every edge intersects a fixed set of $t-1$ vertices.
There has been recent progress on this conjecture; see \cite{AHS12,AFH12,FLM,Fr13,Fr17,HLS,LM}.
In particular,  Frankl \cite{Fr13}
proved that if  $n\geq (2t-1)k-(t-1)$ and $\nu(H)<t$ then $e(H)\le
{n\choose k}-{n-t+1\choose k}$, with further improvement  by Frankl
and Kupavskii \cite{FK18}.

There are  attempts to extend the above conjecture of  Erd\H{o}s to rainbow setting. Let $\mathcal{H} = \{H_1,\ldots, H_t\}$ be a system
of hypergraphs. A matching $\{e_1,\ldots,e_t\}$ is {\it rainbow} if there is a bijection $\phi: [t] \to [t]$ such that $e_i \in E(H_{\phi(i)})$ for each $i \in [t]$.
 In this case, we  say that ${\cal H}$ or
$\{H_1,\ldots, H_t\}$  admits a rainbow matching. Aharoni and
Howard \cite{AH} (see also Huang, Loh and Sudakov \cite{HLS})
conjectured the following rainbow version of the Erd\H{o}s Matching
Conjecture: If $E(H_i)\subseteq {[n]\choose k}$ and $e(H_i)> \max\left\{{n\choose k}-{n-t+1\choose k},{kt-1\choose k}\right\}$ for $i\in [t]$, then
$\{ F_1,\ldots, F_t\}$ admits a rainbow matching.
For the case of 2-graphs, this conjecture was proved by Akiyama and Frankl (see Theorem \ref{AkF}). The case of 3-graphs  was established by Gao, Lu, Ma and Yu \cite{GLMY} provided that n is sufficiently large.
This conjecture was confirmed for $n>3k^2t$ by Huang, Loh and Sudakov \cite{HLS}. The lower bound for  $n$ was improved to $12kt \log(e^2t)$
by Frankl and Kupavskii \cite{FK20}, and subsequently to $3ekt$ by Kupavskii \cite{Ku23}. Lu, Wang and Yu \cite{LWY20} further improved this bound to $n > 2kt$ for $t$ sufficiently large.
Recently, Keevash, Lifshitz, Long and Minzer \cite{KLLM2} proved a general result with $n \ge Ckt$ for some constant $C$ using sharp threshold techniques developed in \cite{KLLM1}.

The notation in this paper is standard. Let $H$ be a $k$-graph. For  $T\subseteq V(H)$, the  neighborhood $N_H(T)$ of $T$  in $H$ is the set $\{A\ |\  A\cup T \in E(H)\}$ and  the \textit{degree}  $d_H(T)$ is $|N_H(T)|$. For any integer  $\ell \ge 0$, the \textit{minimum $\ell$-degree} of $H$ is defined as $\delta_\ell(H):=
\min\{d_H(T): T\in {V(H)\choose \ell}\}$. Note that  $\delta_0(H)=e(H)$ and $\delta_1(H)$ is usually called the \textit{minimum vertex degree} of $H$.
Let $\Delta(H) :=\max \{d_H(v): v\in V(H)\}$ and $\Delta_2(H):=\max \{|N_H(T)|: {T\in \binom{V(H)}{2}}\}$.
For the case of $T=\{u\}$, we simply write $N_H(u)$ and $d_H(u)$ for  $N_H(T)$ and $d_H(T)$, respectively.
To avoid confusion, we sometimes identify  $N_H(T)$ with  the $(k-|T|)$-graph with vertex set $V(H)-T$ and edge set $N_H(T)$.
For any $S \subseteq V(H)$, let $H[S]$ denote the
subgraph of $H$ with $V(H[S])=S$ and $E(H[S])=\{e\in E(H): e \subseteq
S\}$, and let $H-S:= H[V(H)\setminus S]$.  When $T\subseteq U\subseteq V(H)$,  we simply write $d_{U}(T)$ instead of $d_{H[U]}(T)$.



The rest of this paper is organized as follows. In Section 2, we will introduce  necessary notation and definitions. In addition, we will recall a number of lemmas, and  concepts. Section 3 is devoted to the proof of the base case $r=3$. In Section 4, we will analyze the situation for larger matchings. The case where the hypergraph is close to the extremal configuration is discussed in Section 5. Finally, the case far from the extremal configuration is proved in Section 5.

\section{Preliminaries}



Let $e=\{a_1,...,a_k\}$ and $f=\{b_1,...,b_k\}$ be two subsets in $\binom{[n]}{k}$ with $a_1<a_2<...<a_k$ and $b_1<b_2<...<b_k$.
We write $e\le f$ if $a_i\le b_i$ for all $1\le i \le k$, and $e<f$ if $e\le f$ and $e\ne f$.
\begin{myDef}
A $k$-graph $H \subseteq \binom{[n]}{k}$  is {\bf stable}   if $e<f\in F$ implies $e\in H$.
\end{myDef}

Assume $n\geq 2m+2$ and let $\mathcal{G} = \{G_1,...,G_{m+1}\}$ be a system of graphs with the same vertex set $[n]$.
To prove Theorem~\ref{main}, we convert
this rainbow matching problem on graphs  to a matching problem for a special class of $3$-graphs.
We define two special $3$-graphs as following. Write $n-2(m+1) = 2t+\ell $ for integers $t \geq 0$ and $0\le \ell\leq 1 $.
Let $\mathcal{V}= \{v_1,...,v_{m+1}\}$ and $\mathcal{U}= \{u_1,...,u_t\}$ be two disjoint sets such that $(\mathcal{V}\cup \mathcal{U})\cap [n]=\emptyset$. 

\begin{myDef}\label{Def:H(F)}
Let $H(\mathcal{G})$ be the $3$-graph with the vertex set $[n]\cup \mathcal{V}$ and the edge set $\bigcup_{i=1}^{m+1} \{e\cup\{v_i\}\ |\ e\in G_i\}$.
 Additionally, we define $H^\ast(\mathcal{G})$ to be the $3$-graph with the vertex set $[n]\cup \mathcal{V}\cup \mathcal{U}$ and the edge set $E(H(\G)) \cup \bigcup_{i=1}^{t}\{e\cup \{u_i\}\ |\ e\in \binom{[n]}{2}\}$.
\end{myDef}

Note that $\mathcal{G}$ admits a rainbow matching if and only if $H(\mathcal{G})$ has a matching of size $m+1$, which is in turn equivalent to $H^*(\mathcal{G})$ containing a matching of size $m+1+t$.

Let $S(n,m)$ and $D(n,m)$ be two graphs
on the common vertex set $[n]$, defined by
$S(n,m):=G(n,m,r)$ and $D(n,m):=T(2m+1,r)\cup(n-2m-1)K_1$.
 For positive integers $m_1$ and $m_2$,
let $K_{m_1,m_2}$ denote the complete
bipartite graph with part sizes $m_1$ and $m_2$.
If $G_i=K_{m,n-m}$ for every
$1\leq i\leq m+1$, we write
$H(G)=H_S(n,m)$. Similarly, if $G_i=D(n,m)$, then we write $H(\mathcal{G})$ as $H_D(n,m)$.

Let $H_1$ and $H_2$ be two $k$-graphs on the same vertex set $V$ and let $\varepsilon$ be a positive constant.
The following definition gives a quantitative description when $H_2$ contains ``most'' edges of $H_1$.
\begin{myDef}
We say $H_2$ is {\bf $\varepsilon$-close} to $H_1$ if $|E(H_1) \setminus E(H_2)|\le \varepsilon|V|^k$.
\end{myDef}
A \emph{fractional matching} in a $k$-graph $H$ is a function $w : E(H)\to [0, 1]$ such that
\begin{equation}\label{fram}
\sum_{e\in E(H): v\in e} w(e) \le 1
\end{equation}
for any vertex $v \in V(H)$. The fractional matching number $\nu_f(H)=\max \sum_{v \in V(H)} w(v)$ in which the function $w$ satisfies \eqref{fram}.
A fractional matching is called \emph{perfect} if $\sum_{e\in E(H)} w(e) =|V(H)|/k$.

We will frequently use the following lemmas.
\begin{lemma} \label{edgeub}
If $G$ is $K_{r+1}$-free  and the maximum degree of $G$ is $\Delta$, then $e(G) \leq \ex(\Delta,K_{r})+\Delta(n-\Delta)$.
\end{lemma}
\pf Let $v$ be a vertex with degree $\Delta$ and $U=V(G)\setminus N_G(v)$. Note that $e(G)=e(N_G(v))+e(N_G(v),U)+e(U)$. Apparently, $e(N_G(v)) \leq \ex(\Delta,K_r)$ as $G$ is $K_{r+1}$-free. In addition, we have $e(N_G(v),U)+e(U) \leq \sum_{u \in U}d_G(u) \leq \Delta(n-\Delta)$. The lemma follows. \qed

\begin{lemma}\label{combturan}
For positive integers $m_1,\ldots,m_p$ and $r \geq 3$, let $m=\sum_{j=1}^p m_j$. The following holds:
\[
\sum_{j=1}^p \ex(m_j,K_r) \leq \ex(m-p+1,K_r).
\]
\end{lemma}
\pf For each $1 \leq j \leq p$, let $G_j$ be the Tur\'an graph $T(m_j,r-1)$ and fix $v_j \in V(G_j)$. Let $G$ be the graph obtained by identifying
the vertices $v_1,\ldots,v_p$. Note that $G$ has $m-p+1$ vertices and $\sum_{j=1}^p \ex(m_j,K_r)$ edges. Moreover, $G$ is also $K_r$-free. It is clear that $e(G) \leq \ex(m-p+1,K_r)$ and the lemma is proved. \qed

\section{Special Case}

\subsection{Small Matching}
We begin with the following lemma.
\begin{lemma}\label{book-subgraph}
 Assume $n, m, r$ are positive integers with $n>2m$.
Let   $K_m\vee (n-m)K_1$ be the book graph with vertex partition $(A,B)$ such that $|A|=m$ and $|B|=n-m$.
If $G$ is a subgraph of $K_m\vee (n-m)K_1$ and is $K_{r+1}$-free, then $e(G)\leq \ex(m,K_r)+m(n-m)$.
\end{lemma}
\pf Let $d':=\max\{d_G(x)\ |\ x\in B\}$ and  $d:=\max\{\Delta(G[A]),d'\}$. Assume $r \leq m$ for a moment. Since $G$ is $K_{r+1}$-free,  we have
\begin{align*}
e(G)&\leq \ex(d,K_r)+d'(n-d')\\
    &\leq \ex(d,K_r)+d(n-d)\\
&\leq \ex(m,K_r)+m(n-m),
\end{align*}
with equality if and only if $d'=d=m$. Here the last two inequalities follow from the monotonicity of functions.  If $r>m$, then the assertion follows trivially as $\ex(m,K_r)=\binom{m}{2}$.  This completes the proof. \qed

For small matching, we have the following theorem.
\begin{theorem}\label{small-Matching}
Let $n,r,m$ be three positive integers with  $n\geq 10m$ and let   $\G=\{G_1,\ldots,G_{m+1}\}$ be a system of graphs on the same vertex set $V$ with $n$ vertices. If the graph $G_i$ is  $K_{r+1}$-free and
\[
e(G_i)>\ex(n,\{K_{r+1},M_{m+1}\})= \max\{e(G(n,m,r)),e(T(2m+1,r))\}
\]
for each $i \in [m+1]$,  then $\G$ admits a rainbow matching.\end{theorem}

\pf Notice that $\max\{e(G(n,m,r)),e(T(2m+1,r))\}=e(G(n,m,r))$  follows from the assumption $n \geq 10m$.  Let $ A := \{v_1, \ldots, v_{m+1}\} $ and $ V := \{u_1, \ldots, u_n\}$. We define a bipartite graph $ F $ with vertex classes  $A$  and  $V$, where  $ v_i u_j \in E(F) $ if and only if $ d_{G_i}(u_j) > 2m $.

\medskip
\noindent
\textbf{Claim 1.}~For each $i\in [m+1]$ and any $S\subseteq V$ with $|S|\leq m-1$, either $\nu(G_i-S)>2(m+1-|S|)$ or $\Delta(G_i-S)>2m.$

 \medskip
 \noindent
 {\bf Proof of Claim 1:} Suppose the assertion is not true.  Then there exist  $i\in [m+1]$ and a subset  $S\subseteq  V$ with $|S|\leq m-1$ such that $\nu(G_i-S)\leq 2(m+1-|S|)$ and $\Delta(G_i-S) \leq 2m.$ Let $s:=|S|$.  As $G_i$ is $K_{r+1}$-free, Lemma \ref{book-subgraph} implies
\[
e(G_i[S])+e(S,V(G_i)-S)\leq \ex(s,K_r)+s(n-s).
\]
By the assumption $n\geq 10m$ and $e(G_i)>e(G(n,m,r)=\ex(m,K_r)+m(n-m)$, we get
\begin{align}\label{sec3-eq1}
 e(G_i-S)\geq (m-s)(n-s-m)\geq 8m(m-s).
\end{align}
Let $M$ be a maximum matching in $G_i-S$. Then we have
\begin{align}\label{sec3-eq2}
e(G_i-S)\leq \sum_{x\in V(M)}d_{G-S}(x)\leq 4(m+1-s)m.
\end{align}
For $s \leq m-2$, one can check $8m(m-s)> 4(m+1-s)m$ directly. For $s=m-1$, we observe that $(m-s)(n-s-m) > 4(m+1-s)m$. We get a contradiction in each case and this completes the proof of Claim 1. \qed

Let $\nu(F)=s$. Without loss of generality, we may select
$s$  distinct vertices $u_1,\ldots,u_s$ satisfying $d_{G_i}(u_i)>2m$ for each $1 \leq i \leq s$ and write $S:=\{u_1,\ldots,u_s\}$.

If $s=m+1$, then the existence of a rainbow matching $M_{m+1}$ follows from the definition of $F$.
We next consider the case of $s=m$. We claim $e(G_{m+1}-S)>0$. Otherwise, $G_{m+1}$ is a subgraph of $K_m \vee (n-m)K_1$ and Lemma \ref{book-subgraph} implies that $e(G_{m+1}) \leq \ex(m,K_r)+m(n-m)$. This contradicts the assumption $e(G_{m+1})>e(G(n,m,r))$ and the assertion follows.
Now fix an arbitrary edge $e_0\in E(G_{m+1})$.
As $d_{G_i}(u_i)>2m$ for $1 \leq i \leq m$, we can greedily select $m$ pairwise disjoint edges such that $e_i=\{u_i, w_i\}\in E(G_i)$ for each $i$.    Consequently, the edge set $\{e_0,e_1,\ldots,e_m\}$ forms  a rainbow matching of $\{G_1,\ldots,G_{m+1}\}$.


It remains to consider the case where $s<m$. For $s+1\leq j\leq m+1$,  we have
$\Delta(G_{j+1}-S)\leq 2m.$  Otherwise, $u_{j+1}$ is adjacent to $v_{j+1}$ in $F$ and we can find a matching of size $s+1$ in $F$, a contradiction.
 Claim 1 implies that $\nu(G_{j+1}-S)>2(m+1-s)$ for $s+1\leq j\leq m+1$.
Then we may greedily choose $m+1-s$ pairwise disjoint edges $e_{s+1},\ldots,e_{m+1}$ such that $e_j\in E(G_j-S)$  for $s+1\leq j\leq m+1$. Let $R:=\cup_{i=s+1}^{m+1}e_i$.
As in the previous case, we  construct $s$ pairwise disjoint edges $e_1\in E(G_1-R),\ldots,e_s\in E(G_s-R)$ such that $\{e_1,\ldots,e_{m+1}\}$ is a rainbow matching of $\{G_1,\ldots,G_{m+1}\}$. This completes the proof. \qed

\subsection{The case of $r=2$}
We start with recalling the follow theorem by Moon.

\begin{theorem}[Moon, \cite{M}]\label{moon}
    Let $n$ and $m$ be two positive integers such that $n\geq 2m+2$. We have that
\[
\ex(n,\{K_{3},M_{m+1}\})=m(n-m),
\]
where the unique extremal graph is the complete bipartite graph $K_{m,n-m}$.
\end{theorem}

We are now ready to prove Theorem \ref{main}  for the case $r=2$, which we formulate as the following theorem.
\begin{theorem}\label{k=2}
Let $n$ and $m$ be two positive integers such that $n\geq 2m+2$ and $\G=\{G_1,\ldots,G_{m+1}\}$ be a system of graphs on the same vertex set $[n]$.  If   $G_i$ is a $K_3$-free graph and
\[
e(G_i)> \ex(n,\{K_{3},M_{m+1}\})=m(n-m)
\]
for each $i \in  [m+1]$,
then $\G$ admits a rainbow matching.
\end{theorem}
\pf The proof proceeds by induction on  $m$. For the base case in which  $m=1$, Theorem \ref{moon} yields that $\ex(n,\{K_{3},M_2\})=n-1$.
Thus $e(G_i)\geq n$   for $i\in [2]$ by the assumption.   Let $uv\in E(G_1)$. If $G_2-\{u, v\}$ contains an edge $f$, then the set $\{uv, f\}$ forms a rainbow matching. Therefore, we can assume that the set $\{u, v\}$ is a vertex-cover of $G_2$. Since $G_2$ is $K_3$-free, it follows that $N_{G_2}(u)\cap N_{G_2}(v)=\emptyset$ and there are no edges inside $N_{G_2}(u) \cup N_{G_2}(v)$ as  $\{u, v\}$ is a vertex-cover of $G_2$.
Consequently, we have $e(G_2)\leq n - 1$, which leads to a contradiction. This completes the proof of the base case.

  Assume that the assertion holds for   $m\geq 1$. We aim to prove that it also holds for $m'=m+1\leq (n-2)/2$. As we assume $e(G_{m'+1})>m'(n-m')$, it follows that $\Delta(G_{m'+1})\geq m'+1$. We now claim $\Delta(G_{i})\leq n-m'-1$ for each $i \in [m'+1]$. Otherwise,  there exists a vertex $v\in [n]$ and some $i \in [m'+1]$ such that $d_{G_{i}}(v)=\Delta(G_{i})\geq n-m'$. Since $G_{i}$ is $K_3$-free, we have
  \[
  e(G_{i})\leq \Delta(G_{m'+1})(n-\Delta(G_{m'+1}))\leq (n-m')m',
  \]
a contradiction and the claim is true. The last inequality follows from the monotonicity of the function.  For the graph $G_{m'+1}$,  we may assume that $d_{G_{m'+1}}(n)=\Delta(G_{m'+1})$. For each $i\in [m']$, if we let $G_{i}'=G_i-\{n\}$,  then we have
\begin{align*}
e(G_{i}')&\geq e(G_i)-(n-m'-1)\\
&> m'(n-m')-(n-m'-1)\\
&>(n-m')(m'-1)=((n-1)-(m'-1))(m'-1).
\end{align*}
Here the first inequality holds as we already showed that  $\Delta(G_{i})\leq n-m'-1$.
 By induction hypothesis,  $\{G_1',\ldots,G'_{m'}\}$ admits a rainbow matching $M$. As $G_{m'+1}$ is $K_3$-free, we have $e_{G_{m'+1}}(V(M),\{n\})\leq m'$.
 Recall $d_{G_{m'+1}}(n)=\Delta(G_{m'+1})\geq m'+1$. There exists an edge $e\in E(G_{m'+1})$ containing $n$ as an endpoint and $e \cap V(M)=\emptyset. $ Now $M\cup \{e\}$ is a rainbow matching of $\{G_{1},\ldots,G_{m'+1}\}$.
This completes the proof. \qed

\subsection{The case for large $m$}

The following theorem will be used in the proof of Theorem \ref{large s}.

\begin{theorem}[\cite{HYM}]\label{perfect}
For all $\gamma>0$ there is an $n_0$ such that for all $n>n_0$, $n\in 3\mathbb{Z}$, the following holds: Suppose $H$ is a 3-partite 3-uniform hypergraph with  $n$ vertices in each class and $$ \delta_1(H)\geq (\frac{5}{9}+\gamma)n^2.$$ Then $H$ contains a perfect matching.

\end{theorem}

\begin{theorem}\label{large s}
Let $\varepsilon,\alpha,\gamma,\beta$ be four positive constants and let $n,r,m$ be three positive integers such that $r\geq 3$, $ (1/2-\varepsilon)n \leq m \leq (n-2)/2 $ and $1/n\ll \varepsilon\ll \gamma\ll\alpha\ll \beta\ll  1/r$.   If  $\G=\{G_1,\ldots,G_{m+1}\}$ is a system of graphs on the same vertex set $[n]$ such that  $G_i$ is $K_{r+1}$-free and
\[
e(G_i)> \ex(2m+1,\{K_{r+1},M_{m+1}\})
\]
for each $i \in [m+1]$,
then $\G$  admits a rainbow matching.
\end{theorem}

\pf Observe that $\ex(2m+1,\{K_{r+1},M_{m+1}\})=e(T(2m+1,r)$ as we assume $ (1/2-\varepsilon)n \leq m \leq (n-2)/2 $.
Let
 ${\cal V}=\{v_1,...,v_{m+1}\}$ be a set of $m+1$ vertices such that ${\cal V}\cap [n]=\emptyset$. We define $H(\mathcal{G})$ as a 3-graph  whose vertex set is ${\cal V}\cup [n]$ and the edge set is $\bigcup_{i=1}^{m+1} \{e\cup\{v_i\}\ |\ e\in G_i\}$. To prove the desired result, it suffices to show that the
3-graph $H(\mathcal{G})$ has a matching of size $m+1$. Notice that $e(H(\mathcal{G}))=\sum_{i=1}^{m+1} e(G_i)$.

\medskip
\noindent
\textbf{Claim 1.}~For  $i\in [m+1]$, $\Delta(G_i)< (1-1/r+\gamma)n$.

\noindent
{\bf Proof of Claim 1:} Suppose Claim 1 is not true.  Without loss of generality, we assume that $\Delta(G_{1})\geq (1-1/r+\gamma)n$. Let
$v$ in $G_1$ such that $\Delta(G_{1})=d_{G_1}(v)$. Since $G_1[N_{G_1}(v)]$ is $K_{r}$-free, the number of edges in $G_1$  can be bounded as $e(G_1)\leq \ex(d_{G_1}(v),K_{r})+(n-d_{G_1}(v))d_{G_1}(v)$ which can be viewed as the number of edges in a complete $r$-partite graph $K_{n_1,\ldots,n_r}$, where $n_i=d_{G_1}(v)/(r-1)$ for $1 \leq i \leq r-1$ and $n_r=n-d_{G_1}(v)$. One can observe that $e(K_{n/r,\ldots,n/r})-e(K_{n_1,\ldots,n_r})$ is minimized when $d_{G_1}(v)= (1-1/r+\gamma)n$. For  $d_{G_1}(v)= (1-1/r+\gamma)n$, we have
\begin{align*}
e(K_{n/r,\ldots,n/r})-e(K_{n_1,\ldots,n_r})&=n^2 \left(\frac{r-1}{2r}
-\frac{r-2}{2r-2}\left(\frac{1}{r}-\gamma\right)
-\left(\frac{1}{r}-\gamma\right)
\left(1-\frac{1}{r}+\gamma\right)\right)\\
&=n^2\left(\gamma^2
+\frac{(r-2)(3r-2)}{2r(r-1)}\gamma
+\frac{(r-2)(r^2-3r+1)}{2r^2(r-1)}\right)\\
&\geq \frac{r-1}{r}\gamma^2n^2.
\end{align*} Therefore,
\begin{align*}
e(G_1)&\leq \ex(d_{G_1}(v),K_{r})+(n-d_{G_1}(v))d_{G_1}(v)\\
&<\ex(n,K_{r+1})-\frac{r-1}{r}\gamma^2n^2\\
&<\ex(2m+1,K_{r+1})\quad \mbox{(since $\varepsilon\ll \gamma$ and $2m\geq (1-2\varepsilon)n$)},
\end{align*}
a contradiction.  This completes the proof of Claim 1. \qed

Define a set $B$  as follows:
\[
B:=\left\{v\in [n]\ |\ d_{H(\mathcal{G})}(v)<(\frac{5}{9}+\beta){n\choose 2}\right\},
\]
and let $b:=|B|$.

\medskip
\noindent
\textbf{Claim 2.}~$|B|\leq \alpha n$.

\noindent
{\bf Proof of Claim 2:}
Otherwise, assume that  $|B|>\alpha n$. Then we have
\begin{align*}
 e(H(\mathcal{G}))&\leq \frac{1}{2}\left( \sum_{x\in [n]-B}d_{H(\mathcal{G})}(x)+\sum_{y\in B}d_{H(\mathcal{G})}(y)\right)\\
 &<\frac{1}{2}(m+1)(n-\alpha n)(1-1/r+\gamma)n+\frac{1}{2}\alpha n (\frac{5}{9}+\beta){n\choose 2}\quad\mbox{(by Claim 1)}\\
 &<\frac{1}{2}(m+1)n^2(1-\frac{1}{r}+\gamma)-\frac{1}{72}\alpha n^3 \quad \mbox{as $r \geq 3$ and $\varepsilon\ll \alpha\ll \beta\ll 1/r$ }\\
 &<\frac{1}{2}(m+1)(2m+1)^2(1-\frac{1}{r})-\varepsilon n^3\quad\mbox{(since $\varepsilon\ll \alpha$ and $m \geq (1/2-\varepsilon)n)$}\\
 &<(m+1)\ex(2m+1,K_{r+1}),
\end{align*}
a contradiction. This completes the proof of Claim 2. \qed

Without loss of generality, we can assume that the vertices of the hypergraph $H(\mathcal{G})$ are labeled from $1$ to $n$ such that  $d_{H(\mathcal{G})}(1)\geq\cdots\geq d_{H(\mathcal{G})}(n)$.
Set $p:=\lfloor\frac{1}{8\varepsilon}\rfloor$.

\medskip
\noindent
\textbf{Claim 3.}~$d_{H(\mathcal{G})}(2m+2)>0$.

\noindent
{\bf Proof of Claim 3:}
 If $d_{H(\mathcal{G})}(2m+2)=0$, then it follows that   $e(G_i)\leq e(T(2m+1,r))$  for $1\leq i\leq m+1$. This is a contradiction to the assumption for $e(G_i)$,  which proves Claim 3. \qed

\medskip
\noindent
\textbf{Claim 4.}~$d_{H(\mathcal{G})}(2m+1)\geq pn$.

\noindent
{\bf Proof of Claim 4:}
Otherwise, we have
\begin{align*}
  e(H(\mathcal{G}))&\leq (m+1)\ex(2m,K_{r+1})+\sum_{y\in [2m+1,n] }d_{H(\mathcal{G})}(y)\\
  &\leq   (m+1)\ex(2m,K_{r+1})+ p(n-2m)n\\
  & \leq  (m+1)\ex(2m,K_{r+1})+ 2\varepsilon p n^2 \quad\mbox{(since $m \geq (1/2-\varepsilon)n)$}\\
  &<(m+1)\ex(2m+1,K_{r+1}),
\end{align*}
a contradiction. The last inequality is true as $\ex(2m+1,K_{r+1})-\ex(2m,K_{r+1}) \geq 2m-2m/r$, $r \geq 3$, and $\varepsilon$ is small. This completes the proof of Claim 4. \qed

Recall the definition of the set $B$ and notice that $b=|B|$.

\medskip
\noindent
\textbf{Claim 5.}~$H(\mathcal{G})$ has a matching $M$ such that $|M|\leq b$ and $B\cap [2m+2]\subseteq V(M)$.

\noindent
{\bf Proof of Claim 5:} By  the assumption for the degree sequence of $H(\G)$, we may assume $B\cap [2m+2] = [2m+3-s, 2m+2]$, where $s=|B \cap [2m+2]|.$ We  claim that $d_{H(\mathcal{G})}(i)>2(2m+2-i)n\}$ for each $2m+3-s \leq  i \leq 2m+2$. Otherwise, there exists some vertex $i$ with $d_{H(\mathcal{G})}(i)\leq 2(2m+2-i)n$.  Let $j$ be the smallest index with this property. Claims 3 and 4 imply that $j \leq 2m$. Then we have
\begin{align*}
e(H(\mathcal{G}))&\leq  (m+1)\ex(j,K_{r+1})+\sum_{i= j+1}^n d_{H(\mathcal{G})}(i)\\
&\leq  (m+1)\ex(j,K_{r+1})+2(2m+2-j)n(n-j)\\
&< (m+1)\ex(2m+1,K_{r+1})-(m+1)(2m+1-j)\frac{n}{2}+2(2m+2-j)n(n-j)\\
&<(m+1)\ex(2m+1,K_{r+1})-(m+1)(2m+1-j)\frac{n}{2}+3\alpha n^2(2m+2-j) \\
&<(m+1)\ex(2m+1,K_{r+1}),
\end{align*}
which is a contradiction. For the third inequality, to get the Tur\'an graph $T_r(2m+1)$, we add vertices to the Tur\'an graph $T_r(j)$ one by one. Note that each new vertex brings at least $(r-1)n/r \geq n/2$ new edges. We next construct the desired matching greedily.  By Claim 3, there exists an edge $e_0$ containing $2m+2$. Let $M_0=\{e_0\}$. Suppose we have constructed a matching $M_i$ cover  all vertices in $[2m+2-i, 2m+2]$.
If $2m+2-(i+1) \in V(M)$, then let $M_{i+1}=M_i$. For the case  $2m+2-(i+1) \not \in V(M_i)$,
as we proved  $d_{H(\mathcal G)}(2m+2-(i+1))>2(i+1)n$,  we can select an edge $e_{i+1}$ containing $2m+2-(i+1)$ such that $e_{i+1}$ and $M_i$ are disjoint. We set $M_{i+1}=M_i \cup e_{i+1}$ in this case.
 Repeating this process and $M=M_s$ is the desired matching. \qed


Let $m=|M|$ and $X_0:={\cal V}-V(M)$.  By removing vertices with small degree if it is necessary, we define  $Y_0 \subseteq [2m+2]-V(M)$  such that $|Y_0| = 2|X_0| \geq n - 2\alpha n$.
Let  $F:=H(\mathcal{G})[X_0\cup Y_0]$.

Let $R$ be a random subset of $Y_0$ obtained by including each vertex independently with probability $1/2$.
Now consider that 3-partite 3-graph $H_R$ with vertex classes $X_0,R,Y_0-R$ and edge set $E(H_R)=\{e\in E(F)\ |\ |e\cap R|=1\}$.

 For every edge $e\in E(F)$, we define the indicator variable
\[
I_e=\left\{
  \begin{array}{ll}
    1, & \hbox{if $|e\cap R|=1$;} \\
    0, & \hbox{otherwise.}
  \end{array}
\right.
\]
Now for every $x\in X_0\cup Y_0$,
\[
\mathbb{E} (d_{H_R}(x))=\frac{1}{2}d_{F}(x)>(\frac{5}{9}+\beta/2)(n/2)^2.
\]
We apply Janson's inequality (Theorem 8.7.2 in \cite{AS}) to bound the deviation of $d_{H_R}(x)$. For $x\in X_0\cup Y_0$,
         \[\Lambda_x = \sum_{|e\cap f|\geq 2\mbox{ and }x\in e\cap f} \mathbb{P}(I_e=I_f=1)\le \frac{1}{2}n^2(3n/2)\le n^3.
         \]
By Janson's inequality, for any $\beta_0>0$,
       \[
\mathbb{P}(d_{H_R}(x)\le(1-\beta_0)\mathbb{E}(d_{H_R}(x)))\le e^{-\beta_0^2 \mathbb{E}(d_{H_R}(x))/(2+\Lambda/\mathbb{E}(d_{H_R}(x))}=e^{-\Omega(n)}.
         \]
We choose $\beta_0$ such that $\beta_0<\beta/4$. Then we have
\[
\mathbb{P}(d_{H_R}(x)\le (\frac{5}{9}+\beta/4)(n/2)^2\le e^{-\Omega(n)}.
         \]
         The  union bound shows that, with probability
         at least $1-n e^{-\Omega(n)}$, for all $x\in X_0\cup Y_0$,
         $$d_{H_R}(x)\ge (\frac{5}{9}+\beta/4)(n/2)^2.$$
On the other hand,
\[
\mathbb{P}(||R|-|Y_0|/2|\geq n^{3/4})\leq e^{-\Omega(n/2)}.
\]
Thus with positive probability, there exists $R_0\subseteq Y_0$ such that the following two statements hold.
\begin{itemize}
  \item [(i)] $||R_0|-|Y_0|/2|\leq n^{3/4}$;
  \item [(ii)] $d_{H_{R_0}}(x)\geq (5/9+\beta/2)(n/2)^2$.
\end{itemize}

We fix such  $R_0$ such that (i) and (ii) hold.  Without loss generality, assume that $|R_0|\leq |Y_0|/2$. Let $b_1:=|Y_0|-2|R_0|$.  Since $\delta(H(\mathcal{G}))>(5/9){n\choose 2}$, then for any $x\in Y_0-R_0$, there are at least $2b_1n$ edges $f$ such that $x\in f$ and $|f\cap (Y_0-R_0)|\geq 2$. Now we may greedily find a matching $M'$ of size $b_1$ in $H(\mathcal{G})$ such that $V(M')\cap R_0=\emptyset$.
Let $H_2:=H_{R_0}-V(M')$.
Then $\delta(H_2)>(5/9+\gamma)(n/2)^2$. By Theorem \ref{perfect}, $H_2$ has a perfect matching $M''$. It follows that $M \cup M' \cup M''$ is a matching of $H(\mathcal{G})$ of size $m+1$.
This completes the proof. \qed

\section{Close to Extremal Graphs}

\subsection{Close to Space Extremal Graphs}
Let $H_1$ and $H_2$ be two hypergraphs on the same vertex set $V$. For a constant $\alpha>0$, a vertex $v\in V$ is said to be $\alpha$-good for $H_1$ with respect to $H_2$ if
\[
\big|N_{H_2}(v) \setminus N_{H_1}(v)\big| \le \alpha \binom{|V|}{r-1}.
\]
If this inequality fails to hold, we call $v$ an $\alpha$-bad vertex for $H_1$ with respect to $H_2$.

Let $U\cup W$ be a partition of $[n]$ that corresponds to the definition of $H_S(n,m)$, where $|W|=m$ and
 $|U|=n-m$.

\begin{lemma}\label{good-1}
Let $n,m$ be positive integer and let $c<1/10$ be small constant such that $10m\geq n\geq 2m$. Assume $H$ is a  $(1,2)$-partite $3$-graph with vertex partition  $\mathcal{V}\cup [n]$, where $\mathcal{V}:=\{v_1,\ldots,v_{m}\}$. Let $\alpha$ be a constant such that $0<\alpha <10^{-12}$. If every vertex in $H$ is $\alpha$-good with respect to  $H_S(n,m)$,
then $\nu(H)\ge m$.
\end{lemma}

\pf
Let $M$ be a matching in the hypergraph $H$ consisting of edges that each contains exactly one vertex from $\mathcal{V}$, one from $W$, and one from $U$. Among all such matchings, choose $M$ to be of maximum size.
Let $\mathcal{V}':=\mathcal{V}\setminus V(M)$, $U'=U\setminus V(M)$, and
$W'=W\setminus V(M)$.

We claim that $|M|\geq m/2$. Otherwise,
$|\mathcal{V}'|>m/2,\ |W'| > m/2,\ |U'|> n-2m$. Combining $n \leq 10m$ and   the maximality of $M$, we
have, for any $u\in U'$,
\[
|N_{H_S(n,m)}(u)\setminus N_H(u)|\geq |\mathcal{V}'||W'|> n^2/400>\alpha n^2,
\]
a contradiction.

If $|M|< m$, then $\mathcal{V}',U',W'$ are
all nonempty and we are able to choose three distinct vertices
$v\in \mathcal{V}'$, $u\in U'$ and $w\in W'$. Let $e_1$ and $e_2$ be two distinct edges  from $M$. By the maximality of $M$, there is no matching of size three in $H$ that lies within $H[e_1 \cup e_2 \cup \{v, w, u\}]$ and consists only of edges containing exactly one vertex from each of $\mathcal{V}$, $W$, and $U$. Recall that each 3-set $Q$ such that $| Q \cap {\cal V}|=|Q \cap W|=|Q \cap U|=1$  is an edge in $H_S(n,m)$.
Therefore, there exists $f \in E(H_S(n,m))\setminus E(H)$ such that
$f \subseteq e_1\cup e_2\cup \{v,w,u\} $, $|f \cap e_i|=1$ for $i\in [2]$, $|f \cap
\{v,w,u\}|=1$,  and $f$ contains one vertex from each of the sets $\mathcal{V}$, $W$, and $U$.

Note that there are at least ${m/2 \choose 2}$ different choices for $\{e_1,e_2\}$, which
result in distinct choices for $f$. Thus the number of edges in
$E(H_S(n, m))\setminus E(H)$ containing exactly one vertex from
$\{v,w,u\}$ is at least
 $$ {|M|\choose 2}\ge {m/2\choose 2}>n^2/(1000).$$
This implies that for some $x\in \{v,w,u\}$,
$$|N_{H_S(n,m)}(x)\setminus N_H(x)|>n^2/(3000)>\alpha n^2,$$
a contradiction. \qed

 \medskip


Recall that $\ex(n,\{K_{r+1},M_{m+1}\})=f(n,m,r)$.
\begin{lemma}\label{Phk}
	Let $ n,m,r$ be integers and let $\alpha_0,\varepsilon,c$ be constants such that $r\geq 3$,
	$0<1/n\ll \alpha_0\ll \varepsilon\ll 1/r$  and $10m\geq n\geq 2m+\varepsilon n$.  Suppose $H$ is a  $(1,2)$-partite $3$-graph with vertex partition  $\mathcal{V}\cup [n]$, where $\mathcal{V}=\{v_1,...v_{m+1}\}$. If
  $H$ is $\alpha_0$-close to $H_S(n,m)$ and $d_H(x)>f(n,m,r)$ for all $x\in \mathcal{V}$,
	then $\nu(H)\ge m+1$.
\end{lemma}

\pf  Since $H$ is $\alpha_0$-close to $H_S(n,m)$, all but at most $3\sqrt{\alpha_0}\,n$ vertices of $H$ are $\sqrt{\alpha_0}$-good with respect to $H_S(n,m)$.
Let $\mathcal{V}^{\text{bad}}$, $U^{\text{bad}}$, and $W^{\text{bad}}$ denote the sets of $\sqrt{\alpha_0}$-bad vertices in $\mathcal{V}$, $U$, and $W$, respectively. Consequently, we have $|\mathcal{V}^{\text{bad}}| + |U^{\text{bad}}| + |W^{\text{bad}}| \leq 3\sqrt{\alpha_0}\,n$.

Define $\ell:= \max\{|W^{\text{bad}}|+1,|\mathcal{V}^{\text{bad}}|\}$. Choose subsets $W_0 \subseteq W$ and $\mathcal{V}_0 \subseteq \mathcal{V}$ such that $W^{\text{bad}} \subseteq W_0$, $\mathcal{V}^{\text{bad}} \subseteq \mathcal{V}_0$, and $|\mathcal{V}_0| = |W_0|+1 = \ell$.
Set $V_1 := U \cup W_0$ and $W_1 := W \setminus W_0$. Note that $\ell \leq 3\sqrt{\alpha_0}\,n$. Without loss of generality, we write $\mathcal{V}_0 = \{x_1, \dots, x_l\}$.

\medskip
\noindent
\textbf{Claim 1.}~For  $i\in [m+1]$ and any $S\subseteq [n]$ with $m-3\ell\leq |S|\leq m-1$, either $\nu(N_H(x_i)-S)>2(m+1-|S|)$ or $\Delta(N_H(x_i)-S)>2\ell$.

\medskip
\noindent
{\bf Proof of Claim 1:} Let $G:=N_H(x_i)$ and $s:=|S|$. Let $b':=\max\{e_G(x,S)\ |\ x\in V(G)-S\}$ and let $b:=\max\{\Delta(G[S]),b'\}.$

Recall that $G$ is $K_{r+1}$-free. Repeating the proof of Lemma \ref{book-subgraph}, we have
\begin{align*}
e(G[S])+e_{G}(S,[n]-S)&\leq \ex(b,K_r)+b(n-s).
\end{align*}
Suppose that $\nu(G-S)\leq 2(m+1-|S|)$ and $\Delta(G-S)\leq 2\ell$. Let $M$ be a maximum
matching in $G-S$. Then we have
\begin{align*}
e(G-S)\leq \sum_{y\in V(M)}d_{G-S}(y)\leq 4(m+1-s)\ell.
\end{align*}
Thus we may infer that
\begin{align*}
e(G)&\leq e(G-S)+e_G(S,[n]-S)+e(G[S])\\
&\leq \ex(b,K_r)+b(n-s)+4(m+1-s)\ell \\
&\leq \ex(s,K_r)+s(n-s)+4(m+1-s)\ell\\
&\leq \ex(m,K_r)+m(n-m)-(m-s)(n-2m)+4(m+1-s)\ell\\
&\leq \ex(m,K_r)+m(n-m)  \quad \mbox{($\ell \leq 3\sqrt{\alpha_0}n,\alpha_0\ll \varepsilon, n\geq 2m+\varepsilon n$)}\\
&<f(n,m,r).
\end{align*}
This completes the proof of Claim 1. \qed

\medskip
\noindent
\textbf{Claim 2.}~$\{N_H(x_i)-W_1\ |\ i\in [\ell]\}$ has a rainbow matching $M_1$ of size $\ell$.

\medskip
\noindent
{\bf Proof of Claim 2:} For $1\leq i\leq [\ell]$, let $F_i:=N_H(x_i)-W_1$. Let $A:=\{x_1,\ldots,x_\ell\}$ and $B:=[n]-W_1$. We define $L$ as a bipartite graph with vertex partition $(A,B)$
 such that for $x_i\in A$ and $y_j\in B$, $x_iy_j\in E(L)$ if and only if $d_{F_i}(y_j)> 2(\ell+1).$ Let
 $\{x_1y_1,\ldots,x_{\ell_0}y_{\ell_0}\}$ be a maximum matching in $L$ and let $T:=\{y_1,\ldots,y_{\ell_0}\}$, where $\ell_0\leq \ell$. Here note that $d_{F_i-T}(y_j) \leq 2\ell+1$ for $i,j \geq \ell_0+1$. Otherwise, we can get a larger matching in $F$.

We claim that $\{F_j-T:  \ell_0+1\leq j\leq \ell\}$ admits a rainbow matching $M_0$ of size $\ell-\ell_0$. If $\ell_0=\ell-1$, then $|T\cup W_1|=m-1$.
Since $e(N_H(x_i))>f(n,m,r)$, we get $e(N_H(x_\ell)-(T\cup W_1))>0$ and the existence of  a rainbow matching with one edge follows.
Next we may assume that $\ell_0\leq \ell-2$.
 By Claim 1, for $\ell_0+1\leq j\leq \ell,$ we have $\nu(F_j-T)\geq 2(\ell-\ell_0)$.
Let $f_{\ell_0+1}\in N_H(x_{\ell_0+1})$. Suppose that we have constructed $f_{\ell_0+1}\in  F_{\ell_0+1}-T,\ldots, f_i\in F_i-T$. Since $\nu(N_{H}(x_{i+1}))>2(\ell-\ell_0)$, then there exists an edge $f_{i+1}\in E(F_{i+1})$ such that $f_{i+1}\cap (\cup_{j=\ell_0+1}^i f_j)=\emptyset$. Continuing this process, we may find a desired rainbow matching $M_0$ of size $\ell-\ell_0$.

 Since $d_{F_1-T}(y_1)>2\ell$, we may choose an edge $f_1\in F_1-V(M_0)-(T-y_1)$ containing $y_1$ as an endpoint. Suppose that we have constructed $i$ pairwise disjoint edges $f_1,\ldots, f_i$, where $i<\ell_0$ such that
 for $1\leq j\leq i$, $y_i\in f_i$, $|f_i\cap T|=1$, $f_i\cap V(M_0)=\emptyset.$  Since $d_{F_{i+1}}(y_{i+1})>2\ell$ and $(\cup_{j=1}^{i}e_i)\cup V(M_1)\leq 2\ell$,  we are able to find an edge $f_{i+1}$ containing $y_{i+1}$ such that $f_{i+1}\cap T=\{y_{i+1}\}$ and $f_{i+1}\cap ((\cup_{j=1}^{i}f_i)\cup V(M_0))=\emptyset.$ Continuing this process, we can find a rainbow matching of $\{N_H(x_i)-W_1\ |\ i\in [\ell]\}$. This completes the proof of Claim 2. \qed

Let us continue the proof of Lemma \ref{Phk}.
Let $\mathcal{V}':=\mathcal{V}-V(M_1)$, $W':=W-V(M_1)$ and $U'=U-V(M_1)- U^{bad}$. One can see that $|U'|\geq |W'|=|\mathcal{V}'|$. Write $H':=H[\mathcal{V}'\cup W'\cup U']$. Note that each vertex $x\in V(H')$ is $\sqrt{\alpha_0}$-good with respect to $H_S(n,m)$ and $|N_H(x)\setminus N_{H'}(x)| \leq 9\sqrt{\alpha_0}n^{2}$.
Now, for every $x\in V(H')$, we have
\begin{align*}
|N_{H_S(|U'|+|W'|,|W'|)}(x)\setminus N_{H'}(x)|&\leq | N_{H_S(n,m)}(x) \setminus N_{H'}(x)|\\
&\leq 4\sqrt{\alpha_0}n^2 \\
&<\alpha_0^{1/3}|V(H')|^2.
\end{align*}
Thus for every $x\in V(H')$, $x$ is $\alpha_0^{1/3}$-good on $H_S(|U'|+|W'|,|W'|)$. By Lemma \ref{good-1}, $H'$ has a matching $M_2$ of size $|W'|$. It follows that $M_1\cup M_2$ is a matching of size $m+1$ in $H$. This completes the proof. \qed

\begin{lemma}\label{eq-space-lem}
Let $ n,m,r$ be integers and let $\rho,\alpha, c$ be constants such that $r\geq 3$,
	$0<1/n\ll \rho\ll \alpha \ll c\ll 1/r$  and $10m\geq n\geq 2m+cn$.  Suppose $H$ is a  (1,2)-partite $3$-graph with vertex partition  $\mathcal{V}\cup [n]$, where $\mathcal{V}=\{v_1,...v_{m+1}\}$. If the following three conditions hold, then we have $e(H-W)\geq \alpha^4 n^3$  for any $W\in {[n]\choose m}$.
    \begin{itemize}
\item [(i)] For any $x\in \mathcal{V}$, $N_H(x)$ is $K_{r+1}$-free;

        \item [(ii)] $H$ is not $\alpha$-close to $H_S(n,m)$;

        \item[(iii)]  $d_H(x)>f(n,m,r)-\rho n^2$ for all $x\in \mathcal{V}$.
    \end{itemize}
\end{lemma}

\pf  Suppose that there exists $W\in {[n]\choose m}$ such that $e(H-W)\leq \alpha^4 n^3$. We observe that there are at least $m-\alpha^{3/2} n$ vertices $x\in \mathcal{V}$ such that $e(N_H(x)-W)\leq \alpha^{5/2} n^2$.

Since $H$ is not $\alpha$-close to $H_S(n,m)$,  we have \[
|E(H_S(n,m))\setminus E(H)|\geq \alpha |V(H)|^3=\alpha (n+m)^3.
\]
Thus there are at least $\alpha n/2$ vertices $x\in \mathcal{V}$ such that
$|E(K_{W,[n]-W})\setminus E(N_H(x))|\geq \alpha n^2$. (Note that the partition is $(W,[n]-W)$ with $|W|=m$). It follows that there exists $v\in {\cal V}$ such that
$|E(K_{W,[n]-W}) \setminus N_H(v)|\geq \alpha n^2$
and $e(N_H(v)-W)\leq \alpha^{5/2} n^2$.  Let $G:=N_H(v)$. Thus we may infer that
\begin{align*}
e(G[W])&\geq f(n,m,r)-\rho n^2-m(n-m)+|E(K_{W,[n]-W})\setminus E(N_H(x))|-e(G-W)\\
&\geq \ex(m,K_{r})+\alpha n^2/2,
\end{align*}
which implies that there exists $y\in W$ such that $d_{G[W]}(y)\geq (1-1/(r-1))m+\alpha n$.

Let $b:=\max\{e_G(\{x\},W)\ |\ x\in [n]-W\}$. If $b\leq m-\alpha^2 n$, then we
have
\begin{align*}
e(G)&= e(G[W])+e_G(W,[n]-W)+e(G-W)\\
&\leq (\ex(b,K_{r})+(m-b)m)+b(n-m)+\alpha^{5/2} n^2\\
&\leq (1-1/(r-1))b^2/2+b(n-2m)+m^2+\alpha^{5/2} n^2\\
&< (1-1/(r-1))(m-\alpha^2n)^2/2+(m-\alpha^2n)(n-2m)+m^2+\alpha^{5/2} n^2\\
&<(1-1/(r-1))(m-\alpha^2n)^2/2+m(n-m)-\alpha ^2n(n-2m)+\alpha^{5/2} n^2\\
&<(1-1/(r-1))m^2/2+m(n-m)-c\alpha^2n^2 \quad \mbox{($n \geq 2m+cn$ and $\alpha \ll c$)}\\
&<\ex(m,r)+m(n-m)-\rho n^2 \quad \mbox{ ($\rho \ll \alpha \ll c$)}\\
&\leq f(n,m,r)-\rho n^2,
\end{align*}
a contradiction.

Next we may assume that $b>m-\alpha^2n$. Recall $|E(K_{W,[n]-W}) \setminus N_H(v)|\geq \alpha n^2$ which implies that $e_G(W,[n]-W) \leq m(n-m)-\alpha n^2$. We have
\begin{align*}
e(G)&= e(G[W])+e_G(W,[n]-W)+e(G-W)\\
&\leq \ex(b,K_{r})+(m-b)m+(m(n-m)-\alpha n^2)+\alpha^{5/2} n^2\\
&\leq \ex(m,K_r)+\alpha^2 n^2+(m(n-m)-\alpha n^2)+\alpha^{5/2} n^2\\
&\leq \ex(m,K_r)+m(n-m)-\rho n^2, \quad \mbox{($\rho \ll \alpha \ll 1/r$)}\\
&\leq f(n,m,r)-\rho n^2,
\end{align*}
a contradiction. This completes the proof. \qed

\subsection{Close to Diversity Extremal Graphs}

\begin{theorem}[Nikiforov,\cite{Nik}]\label{stabilityclique}
Let $n,r$ be positive integers and let  $\beta$ be a constant such that $1/n\ll \beta\ll 1/r$.
Let $G$ be a $K_{r+1}$-free graph with $n$ vertices and $(1-1/r-\beta){n\choose 2}$ edges. Then $|E(T(n,r))-E(G)|\leq \beta^{1/4}n^2$.
\end{theorem}

\begin{lemma}\label{Divesity-close}
Let $n,m,r$ be positive integers  and let
$\alpha, \beta, \gamma,c $, such that $r\geq 3,$
	 $0< \alpha\ll \beta\ll \gamma\ll c\ll 1/r
    $,  and $2m+cn \leq  n \leq 10m$.  Let $H$ be a  (1,2)-partite $3$-graph with vertex partition  $\mathcal{V}\cup [n]$, where $|\mathcal{V}|=m+1$. Suppose that for any $x\in \mathcal{V},$ $N_H(x)$ is $K_{r+1}$-free.  If there exists $U\in {[n]\choose 2m+1}$ such that $e_H(\mathcal{V},U,[n]-U)+e(H[{\cal V} \cup ([n]-U]))\leq \alpha n^3$, and $|N_{H}(x)|>f(n,m,r)$ for all $x\in \mathcal{V}$,
	then $\nu(H)\ge m+1$.
\end{lemma}


\pf Fix a set $U_1\in {[n]\choose 2m+1}$ such that
\begin{align}\label{clo-div-eq1}
   e_H(\mathcal{V},U_1,[n]-U_1)+e(H[{\cal V} \cup ([n]-U_1)])\leq \alpha n^3.
\end{align}
 Let $U_2:=[n]-U_1$ and $H_1:=H[U_1\cup \mathcal{V}]$.

Define
\[
\mathcal{V}_1:=\left\{x\in \mathcal{V}\ \left|N_H(x) \setminus \binom{U_1}{2}\ \right|\geq \sqrt{\alpha} n^2\right\} \quad \mbox{ and } \quad  \mathcal{V}_2:=\mathcal{V}-\mathcal{V}_1.
\]
 The inequality (\ref{clo-div-eq1}) implies that  $|\mathcal{V}_1| \leq \sqrt{\alpha} n$. Moreover, for every vertex $x\in \mathcal{V}_2$, we have  $|N_H(x) \cap \binom{U_1}{2}| \geq f(n,m,r)-\sqrt{\alpha}n^2$. Applying Theorem \ref{stabilityclique}, we conclude that
$N_H(x)$ is $\alpha^{1/2}$-close to $T(2m+1,r)\cup (n-2m-1)K_1$.


Recall $H_1=H[U_1\cup \mathcal{V}]$. We define
\[
B:=\{x\in U_1\ |\ d_{H_1}(x)<(1-1/r-2\beta )2m^2\}.
\]
Let $b:=|B|$ and $W_1=U_1\setminus B$. By (\ref{clo-div-eq1}), we
have 
\begin{align}\label{H1-lower}
e(H_1)\geq mf(n,m,r)-\alpha n^3.
\end{align}
Other hand, since $N_H(x)$ is $K_{r+1}$-free for all $x\in \mathcal{V}$, then we have
\begin{align*}
e(H_1)&\leq e(H_1-B)+\sum_{x\in B}d_{H_1}(x)\\
&\leq m(1-1/r){2m+1-b\choose 2}+b(1-1/r-2\beta)2m^2\\
&\leq m(1-1/r)(2m-b)^2/2+b(1-1/r-\beta)2m^2,
\end{align*}
i.e.
\begin{align}\label{H1-upper}
e(H_1)\leq m(1-1/r)(2m-b)^2/2+b(1-1/r-\beta)2m^2.
\end{align}
Combining (\ref{H1-lower}) and (\ref{H1-upper}), we may get $b \leq \sqrt{\alpha}n$.  Let $H_2$ be the subhypergraph of $H$ obtained by deleting the edges in ${\cal V} \cup W_1$.

\medskip
\noindent
{\bf Claim 1.}~$H_2$ has a matching $M$  such that $|V(M)\cap ([n]-W_1)|\geq b+1$.

 Suppose for contradiction that the assertion is not true.
Since $d_{H}(x)>f(n,m,r)$ for every $x\in \mathcal{V}$, we may assume $b\geq 1$.
Let $M$ be a matching in $H_2$ that maximizes $|V(M)\cap ([n]\setminus W_1)|$. Among all such matchings, choose $M$ to minimize $|V(M)\cap W_1|$.
Define $H_3:=H-(\mathcal{V}\cap V(M))$. Let $l:=|M|$.

By the definition of $H_2$, every edge $e\in M$ satisfies $|e\cap W_1|\leq 1$.
Observe that for any $u_1\in \mathcal{V}\setminus V(M)$ and $u_2\in [n]\setminus(W_1\cup V(M))$,
\begin{align}\label{u1u2-empty}
N_{H_2-V(M)}(\{u_1,u_2\})=\emptyset.
\end{align}

Let
\[
M_1:=\{e\in M\mid e\cap W_1=\emptyset\} \qquad \mbox{ and } \qquad M_2:=M\setminus M_1.
\]
By the choice of $M$, for any $x\in V(M_2)\cap ([n]\setminus W_1)$, we have
\begin{align}\label{C1-eq1}
N_{H_2-\big((V(M)\setminus\{x\})\cup W_1\big)}(x)=\emptyset.
\end{align}

We now claim that for any $f=\{v,u_1,u_2\}\in M_1$ with $v\in\mathcal{V}$,
\begin{align}\label{degree-sum1}
d_{H_3}(u_1)+d_{H_3}(u_2)\leq 4\left(1-\frac1r-\beta\right)m(m-l).
\end{align}
Without loss of generality, we may assume
\[
\nu\big(N_{H_2-(W_1\cup V(M-f))}(u_1)\big)\geq \nu\big(N_{H_2-(W_1\cup V(M-f))}(u_2)\big).
\]
Consider that $\nu\big(N_{H_2-(W_1\cup V(M-f))}(u_1)\big)\geq 3$. Then $N_{H_3-V(M-f)-u_1}(u_2)=\emptyset$. Otherwise, suppose that there exists $e_2\in N_{H_3-V(M-f)-u_1}(u_2)$.
We may select $e_1\in N_{H_2-(W_1\cup V(M-f))}(u_1)$ such that $e_1\cap e_2=\emptyset$.
Then
\[
M':=(M\setminus\{f\})\cup \big\{e_1\cup\{u_1\},\,e_2\cup\{u_2\}\big\}
\]
is a matching satisfying
\[
\big|V(M')\cap ([n]\setminus W_1)\big|>\big|V(M)\cap ([n]\setminus W_1)\big|,
\]
contradicting the choice of $M$.
Thus we have
\begin{align*}
d_{H_3}(u_1)+d_{H_3}(u_2)
&< 2\left(1-\frac1r-2\beta\right)m^2 + m(n-2m)+d_{H_3}(u_2)\\
&\leq 2\left(1-\frac1r-2\beta\right)m^2 + m(n-2m)+2|M|m\\
&< 4\left(1-\frac1r-3\beta\right)m(m-l),
\end{align*}
where the last inequality follows from assumptions $2m\leq n\leq 5m/2$ and $b\leq \beta n$.

Hence we may suppose
\[
2\geq \nu\big(N_{H_2-(W_1\cup V(M-f))}(u_1)\big)\geq \nu\big(N_{H_2-(W_1\cup V(M-f))}(u_2)\big).
\]
Again by Hall's Theorem,
\begin{align*}
d_{H_3}(u_1)+d_{H_3}(u_2)
&< 4\left(1-\frac1r-2\beta\right)m^2
+2\max\big\{|\mathcal{V}\setminus V(M)|,\, n-2m+b\big\}\\
&<4\left(1-\frac1r-\beta\right)m(m-l).
\end{align*}

Define
\[
M_{21}:=\big\{g\in M_2\,\big|\, d_{H_2-V(M-g)}\big(g\cap ([n]\setminus W_1)\big)\leq 2n\big\},
\qquad
M_{22}:=M_2\setminus M_{21}.
\]

Take $g=\{v,v_1,v_2\}\in M_2$ with $v\in\mathcal{V}$, $v_1\in W_1$, and $v_2\in [n]\setminus W_1$.
Firstly, consider  $g\in M_{21}$. Recall that $e_{H_2-V(M-g)}\big(v_2,[n]\setminus(W_1\cup V(M)),\mathcal{V}-V(M)\big)=0$.  Consequently,
\begin{align}\label{v1v2-sum}
d_{H_3}(v_1)+d_{H_3}(v_2)
<nm+ 2n+bm
<4\left(1-\frac1r-\beta\right)m(m-l).
\end{align}
Next consider $g\in M_{22}$. By the choice of $M$,
\begin{align}
e_{H_2-V(M-g)}\big(\mathcal{V}\setminus V(M),\,v_1,\,[n]\setminus(V(M)\cup W_1)\big)=0.
\end{align}

Let $\mathcal{V}':=\mathcal{V}\setminus V(M)$.
Set $\ell:=|M|$ and $s:=|M_1\cup M_{21}|$.
We obtain
\begin{align*}
e(H_3)
&\leq \sum_{x\in V(M_1\cup M_{21})}d_{H_3}(x)
+e\Big(H_3\big[\big((W_1\cup V(M_{22}))\setminus V(M_{21})\big)\cup \mathcal{V}'\big]\Big)\\
&\quad +e_{H_3}\big(\mathcal{V}',\, V(M_{22})\cap [n],\, [n]\setminus(V(M)\cup W_1)\big)\\
&=\sum_{x\in V(M_1\cup M_{21})}d_{H_3}(x)
+e\Big(H_3\big[\big((W_1\cup V(M_{22}))\setminus V(M_{21})\big)\cup \mathcal{V}'\big]\Big)\\
&\leq 4s\left(1-\frac1r-\beta\right)m(m-l)
+(m-l)\,e\big(T(2m+1-2s,r+1)\big)\\
&<(m-\ell)\,f(n,m,r),
\end{align*}
which yields a contradiction.
This completes the proof of Claim~1.

By Claim 1, we may choose a matching $\mathcal{M}_1$ such that $|V(\mathcal{M}_1)\cap ([n]-W_1)|\in\{ b+1,b+2\}$. Write $m_1:=|\mathcal{M}_1|$. Let $\mathcal{V}_1':=\mathcal{V}_1-V(\mathcal{M}_1)$ and  $m_2:=|\mathcal{V}_1'|$.
 As the definition of ${\cal V}_1$, we are able to greedily find a matching $\mathcal{M}_2$ of size $m_2$ in $H-V(\mathcal{M}_1)$  such that $\mathcal{V}_1'\subseteq V(\mathcal{M}_2)$.

Note that ${\cal V}_1 \cup {\cal V}_2$ contains at most $2m_1+2m_2-b-1$ vertices in $W_1$.
Thus we can choose $U_1'\subseteq U_1-B-V({\cal M}_1)-V({\cal M}_2))$ such that $|U_1'|=2m+2-2(m_1+m_2)$.
Let $\mathcal{V}_4:=\mathcal{V}-V({\cal M}_1)-V({\cal M}_2)$ and $H_4:=H[\mathcal{V}_4\cup U_1']$. Notice that $|{\cal V}_4|=m+1-(m_1+m_2)$.
Then $d_{H_4}(x)\geq (1-1/r-\gamma)2m^2$  for any $x\in U_1'$.
Observe that a perfect matching  in $H_4$ together with $\mathcal{M}_1 \cup \mathcal{M}_2$ is a matching of size $m+1$ in $H$.

Next we show that $H_4$ contains a perfect matching. Let $ \mathcal{R}_1$ be a random subset of $ U_1' $ where each vertex $ y \in U_1' $ is selected independently with probability $ 1/2 $. Let $\mathcal{R}_2:=U_1'-\mathcal{R}_1$. We thus obtain a random 3-partite 3-graph $ F $ with vertex partition $ (\mathcal{V}_4, \mathcal{R}_1, \mathcal{R}_2) $, where an edge $ e \in E(F) $ if and only if $ e \in E(H_4) $ and $ |e \cap \mathcal{R}_1| = 1 $.

Note that $ \mathbb{E}[d_F(x)] \geq (1 - 1/r - 3\gamma)m^2 $ for all $ x \in V(F) $. By standard probabilistic arguments, $ F $ satisfies the following properties with probability $ 1 - o(1) $:
\begin{itemize}
    \item $ \delta(F) > (1 - 1/r - 4\gamma)m^2 $,
\item  $ (1 - \alpha)m \leq |\mathcal{R}_1| \leq (1 + \alpha)m $,
\item $ e(F - \mathcal{R}_i) \geq m^2/8 $ for each $ i \in [2] $.
\end{itemize}

Fix such a subset $ \mathcal{R}_1 $ and denote the resulting graph by $ H_5 $. Without loss of generality, assume $ |\mathcal{R}_1| \leq |\mathcal{R}_2|$. We can then greedily construct a matching $ \mathcal{M}_3 $ in $ H_5 - \mathcal{R}_1 $ of size $ (|\mathcal{R}_2| - |\mathcal{R}_1|)/2 $.

Finally, observe that $ \delta(H_5 - V(\mathcal{M}_3)) > (1 - 1/r - 6\gamma)m^2 > (5/9 + \alpha)m^2 $. By Theorem \ref{perfect}, $ H_5 - V(\mathcal{M}_3) $ admits a perfect matching $ \mathcal{M}_4 $. It follows that $ \bigcup_{i=1}^4 \mathcal{M}_i $ is a matching of size $ m + 1 $ in $ H $, completing the proof.\qed

\section{Non-close extremal graphs}

\subsection{Absorbing Lemma}

Our strategy to prove Theorem~\ref{main} is to find a small matching $M'$ in $H$ that can be used to ``absorb" small sets of
vertices,  find an almost perfect matching $M''$ in $H-V(M')$, and  then use $M'$ to absorb $V(H)\setminus V(M'\cup M'')$.
In this section, we prove such an  absorbing lemma for $(1,2)$-partite $3$-graphs. Our proof follows along the same lines
as in \cite{RRS06}.

\begin{lemma}\label{absorbing}
 Let $\rho, \rho',\alpha$ be constants such that $0<1/n\ll \rho'\ll \rho\ll \alpha\ll 1$.
  Let $H$ be a $(1,2)$-partite $3$-graph with partition classes $Q,P$ such that $2|Q|=|P|$. If for any $v\in Q$, $d_H(v)>\alpha n^2$ and $Q$ contains at least $\alpha n$ vertices $v$ such that $d_H(v)={n\choose 2}$.
Then $H$ has  a matching $M$ such that  $|M|\leq \rho n$ and,
for any subset $S\subseteq V(H)\setminus V(M)$
with $|S|\leq \rho' n$ and $2|S\cap Q|=|S\cap P|$,
$H[S\cup V(M)]$ has a perfect matching.
\end{lemma}
\pf
We call a balanced 6-element set $A\subseteq V(H)$ an {\it absorbing set} for a balanced $3$-element set $T\subseteq V(H)$
if $H[A]$ has a matching of size 2 and $H[A\cup T]$ has a matching of size $3$.
Denote by $\mathcal{L}(T)$ the collection of all absorbing sets for $T$.

\medskip
\textbf{Claim 1.~}For every balanced  $T\in  {V(H)\choose 3}$, $|\mathcal{L}(T)|> \alpha^4 n^6$.
\medskip

Let  $T=\{u_0,u_1,u_2\}\in  {V(H)\choose 3}$ be balanced, with
$u_0\in Q$ and $u_1,u_2\in P$. We form an absorbing set for $T$ by choosing three pairwise disjoint 2-sets $U_0,U_1,U_2$ in order.

First, we choose a 2-set $U_0\subseteq P\setminus T$ such that $U_0\cup \{u_0\}\in E(H)$. The number of choices
for $U_0$ is  at least
$$d_H(u_0)-2{n-1\choose 1}> \alpha n^2-2n>\alpha n^2/2.$$

Now fix a choice of $U_0$, and let $U_0=\{w_1,w_2\}$. Note that,
for $i\in [2]$,
$$
|N_{H}(u_i)\cap N_{H}(w_i)|\geq (2\alpha n^2)/3 $$
because there are $\alpha n$ vertices $v \in Q$ such that $d_H(v)=\binom{n}{2}$.
For $i\in [2]$, we choose 2-sets $U_i$ from $(V(H)\setminus T)\setminus \bigcup_{j=0}^{i-1}U_j$ such that $U_i\cup \{u_i\}$ and $U_i\cup \{w_i\}$
are both edges of $H$. For each choice of $U_j$, $0\le j\le i-1$, the number of choices for $U_i$ is at least
$$ |N_{H}(u_i)\cap N_{H}(w_i)|-4n \ge \frac{1}{2}\alpha n^2.$$

Let $A=\bigcup_{i=0}^2 U_2$. Then
$\{U_i\cup \{w_i\}: i\in [2]\}$ is a matching in $H[A]$, and $\{U_i\cup \{u_i\}: i\in [2]\cup \{0\}\}$  is a matching in $H[A\cup T]$.
Thus $A$ is an absorbing set for $T$.
Since there are more than $\alpha^3n^{6}/8$ choices of $(U_0,U_1,U_2)$,
there are more than $\alpha^3n^{6}/(16\times 24)$ absorbing sets for $T$.  This completes the proof of Claim 1.

Now, form a family $\mathcal{F}$ of subsets of $V(H)$ by selecting each of the ${n/2\choose 2}{n\choose 4}$
possible balanced 6-sets independently with probability
\[
p=\frac{\rho n}{2{n/2\choose 2}{n\choose 4}}.
\]
Then, it follows from  Chernoff's bound that, with probability $1-o(1)$ (as $n\rightarrow \infty$),
\begin{itemize}
\item [(i)] $|\mathcal{F}|\leq \rho  n$, and
\item [(ii)] $|\mathcal{L}(T)\cap \mathcal{F}|\geq p |\mathcal{L}(T)|/2\geq \alpha^4\rho n/4$ for all balanced $T\in {V(H)\choose 3}$.
\end{itemize}

Furthermore, the expected number of intersecting pairs of sets in $\mathcal{F}$ is at most
\[
{n/2\choose 2}{n\choose 4} \left[2{n/2-1\choose 1}{n\choose 4}+4{n-1\choose 3}{n/2\choose 2} \right]p^2< \rho^{1.5}n.
\]
Thus, using Markov's inequality, we derive that, with probability at least 1/2,
\begin{itemize}
\item [(iii)] $\mathcal{F}$ contains at most $2\rho^{1.5}n$ intersecting pairs.
\end{itemize}

Hence, with positive probability, $\mathcal{F}$ satisfies (i), (ii), and (iii).
Let $\mathcal{F}'$ be obtained from $\mathcal{F}$ by removing one set from each intersecting pair and deleting all non-absorbing sets.
Then  $\mathcal{F}'$ consists of pairwise disjoint absorbing sets, such that, for each $T\in {V(H)\choose 4}$,
\[
 |\mathcal{L}(T)\cap \mathcal{F}'|\geq \alpha^4\rho n/2.
\]

 Since $\mathcal{F}'$ consists only of pairwise disjoint absorbing sets, $H[V(\mathcal{F}')]$ has a perfect matching, say  $M'$. Then $|M'| \le \rho n$.
 To complete the proof, take an arbitrary $S\subseteq V(H)\setminus V(M')$ with $|S|\leq \rho' n$ and $2|S\cap Q|=|S\cap P|$,
 where $\rho'\le \alpha^4\rho$.
Note that $S$ can be partitioned into $t$ balanced $3$-sets, say $T_1,\ldots, T_t$, for some $t\le \rho' n/3<\alpha^4\rho n/2$.
We can greedily
choose distinct absorbing sets $A_i\in \mathcal{F}'$ in order for $i=1, \ldots, t$, such that $H[A_i\cup T_i]$ has a perfect matching.
Hence, $H[S\cup V(M')]$ has a perfect matching as required. \qed

\subsection{Stability Lemma}

\begin{theorem}[Balogh et al., \cite{BBCL}]\label{Kr-free-close}
Let $H$ be a $K_{r+1}$-free graph on $n$ vertices.  If
$e(H)\geq ex(n,K_{r+1})-\rho n^2,$ then $H$ is $\sqrt{\rho}$-close to $T(n,r)$.
\end{theorem}

For the perfect fractional matching part, we need a
result of Berge \cite{Be58} on maximum matchings.
 For a graph $G$, we use $c_o(G)$ to denote the number of odd components in $G$.

\begin{theorem}[Berge] \label{Berge}
Let $G$ be a graph on $n$ vertices. Then there exists $S\subseteq V(G)$ such that every component of $G-S$ is factor-critical  and
$$n-2\nu(G)=odd(G-S)-|S|.$$
\end{theorem}

\begin{lemma}\label{stable-lem}
Let $ \rho,\beta, c$ be constants such that $0<\rho\ll \beta\ll c$, and let $m, n,r$ be positive integers such that $r\geq 3$,
 $cn\leq m\leq (1-\beta)n/2$ and $1/n\ll 1/r$. Let $G$ be a $K_{r+1}$-free 2-graph with vertex set  $V[n]$. If  $\nu(G)\leq m$ and  $e(G)> f(n,m,r)-\rho n^2$, then $G$ is $\rho^{1/4}$-close to $G(n,m,r)$ or $T(2m+1,r)\cup(n-2m-1)K_1$.
\end{lemma}
\pf Since $\nu (G)\leq m$, by Theorem \ref{Berge},
there exists $S\subseteq V(G)$ such that  every component of $G-S$ is factor-critical and
\begin{equation} \label{odd}
n-2\nu(G)= odd(G-S)-|S|.
\end{equation}
Let $q:=odd(G-S)$ and  $C_1,\ldots,C_q$ be the odd components of $G-S$. Thus $\nu(G)=s+\tfrac{1}{2}\sum_{i=1}^q(|C_i|-1) \leq m$ and consequently
\begin{equation}\label{para}
s \leq m, \qquad q \geq n-2m+s, \qquad \sum_{i=1}^q(|C_i|-1) \leq 2(m-s).
\end{equation}
We set $D:=\cup_{i\in [q]}V(C_i)$
   and choose $v_1\in D$ such that $d_G(v_1)\geq d_G(x)$ for all $x\in D$. Without loss of generality, suppose that $v_1 \in V(C_1)$. For $2\leq i\leq q$, we choose $v_i\in V(C_i)$ and let $T:=\{v_2,\ldots,v_q\}$.
Let $F(n,s)$ be a graph with vertex set $[n]$ and edge set
\[
E(G-(D-V(C_1))\cup (\bigcup_{i=2}^q\{v_iu\ |\ u\in N_G(v_i)\cap S\})\cup (\bigcup_{x\in D-(V(C_1)\cup T)} \{xu\ |\ u\in N_G(v_1)\}).
\]
Namely, for each $2 \leq i \leq q$, we keep only  neighbors of $v_i$ in $S$. For each $v \in \cup_{i=2}^q C_i \setminus T$,
we create a copy of $v_1$.
Note that $e(F(n,s))\geq e(G)$, $\nu(F(n,s))\leq \nu(G)$ and $F(n,s)$ is $K_{r+1}$-free.
Moreover, one can see that $q':=odd(F(n,s)-S)=q$. Let $C_1'\ldots,C_q'$ be the odd components of $F(n,s)-S$ such that
$|V(C_1')|\geq \ldots\geq |V(C_q')|$.  Let
\[
b:=\min \{|W|\ |\ W\subseteq S\mbox{ and $G[S]-W$ is $K_r$-free}\}.
\]

\medskip
\noindent
\textbf{Claim 1.~}$0\leq s\leq 8\rho n$ or $s\geq m-8\rho n$.
\medskip

\noindent
{\bf Proof of Claim 1:} Suppose that $8\rho n\leq s\leq  m-8\rho n$.
Since $G$ is $K_{r+1}$-free, for $x\in D-V(C_1')$.
The assumption $\nu(G) \leq m$ together with the equation \eqref{odd} imply that $q \geq n-2m+s$. Thus we have
\begin{align*}
e(G)&\leq e(F(n,s))\\
&\leq (q-1)(s-b)+((n-q-\ell+1)\ell+\ex(\ell,K_{r}))\\
&\leq (n-2m+s-1)(s-b)+((2m-s+1-\ell)\ell+\ex(\ell,K_r))\\
&\leq (n-2m-1+s)(s-b)+(2m-s-\ell+1)\ell+(1-\frac{1}{r-1})\frac{\ell^2}{2}\\
&\leq (n-2m+s)(s-b)+(2m-s)\ell-(1+\frac{1}{r-1})\ell^2/2+\rho^4 n^2.
\end{align*}
 If $(2m-s)(1-1/r)<s-b$, by the monotonicity of the function $f(\ell)=(2m-s)\ell-(1+\frac{1}{r-1})\ell^2/2$, then we get
\begin{align*}
e(G)&\leq (n-2m+s)(s-b)+(2m-s)(s-b)-(1+\frac{1}{r-1})(s-b)^2/2+\rho^4n^2  \\
&=n(s-b)-(1+\frac{1}{r-1})(s-b)^2/2+\rho^4n^2  \\
&\leq ns-(1+\frac{1}{r-1})s^2/2+\rho^4n^2 \\
&<f(n,m,r)-\rho n^2.
\end{align*}
Here the second-to-last inequality can be verified directly by noting $r \geq 3$ and $s \leq m $ and the last inequality follows as $s \leq m-8\rho n$.

If  $(2m-s)(1-1/r)\geq s-b$, then  the maximum value of the function $f(\ell)=(2m-s)\ell-(1+\frac{1}{r-1})\ell^2/2$ is attained at $\ell=(2m-s)(1-1/r)$ and we have
\begin{align*}
e(G)&\leq (n-2m+s)(s-b)+(2m-s)^2(1-1/r)-(1+\frac{1}{r-1})((2m-s)(1-1/r))^2/2+\rho^4n^2\\
&=(n-2m+s)s+(2m-s)^2(1-1/r)/2+\rho^4n^2\\
&<f(n,m,r)-\rho n^2.
\end{align*}
The last inequality holds because either $s< 8\rho n$ or $s \geq m-8\rho n$ as well as the monotonicity of the function.
This is  a contradiction to the assumption again and  the proof of Claim 1 is complete.

Assume  $|V(C_1)|\geq \ldots\geq |V(C_q)|$. Recall that
\[
\nu(G):=s+\sum_{i=1}^q (|V(C_i)|-1)/2.
\]
We have two cases depending on $s$.

\medskip
\noindent
\textbf{Case 1.~}$0\leq s\leq 8\rho n$.
\medskip

Recall $m\geq cn$ and $\rho\ll c$. Note that $\sum_{i=1}^q (|V(C_i)|-1) \leq 2m-2s$ by \eqref{para}. As we aim to show an upper bound for $\sum_{i=1}^q \ex(|C_i|,K_{r+1})$, we may assume  $\sum_{i=1}^q (|V(C_i)|-1) = 2m-2s$.  We  claim $|V(C_1)|\geq 2m-2s-\sqrt{\rho} n$. Otherwise, for each  $2 \leq i \leq q$, we can choose nonempty subset $C_i' \subseteq C_i$ such that $\sum_{i=2}^q |C_i'|=\sqrt{\rho} n+q+1$.
Then we have
\begin{align*}
e(G)&\leq sn+\sum_{i=1}^q \ex(|V(C_i)|, K_{r+1})\\
    & \leq sn+\ex(2m-2s-\sqrt{\rho}  n-1,K_{r+1})+\sum_{i=2}^q \ex(|C_i'|,K_{r+1})\\
    & \leq sn+\ex(2m-2s-\sqrt{\rho}  n-1,K_{r+1})+\sum_{i=2}^q \ex(|C_i'|-1,K_{r+1})\\
    & = sn+\ex(2m-2s-\sqrt{\rho}  n-1,K_{r+1})+\ex(\sqrt{\rho} n+3,K_{r+1}) \\
    &\leq sn+\ex(2m-2s-\sqrt{\rho}  n,K_{r+1})+\ex(\sqrt{\rho} n+2) \\
&\leq sn+(1-1/r)(2m-2s-\sqrt{\rho}n)^2/2+ (1-1/r)(\sqrt{\rho} n+2)^2/2\\
&<\ex(2m+1,K_{r+1})-\rho n^2.
\end{align*}
Here the second and the fifth inequalities follow from the convexity of the function $\ex(n,K_{r+1})
$. Lemma \ref{combturan} yields the third inequality.
We get a contradiction and $|V(C_1)|\geq 2m-2s-\sqrt{\rho} n$ follows.
Now,
\begin{align*}
 e(G[S\cup V(C_1)])&\geq e(G)-e_G(S,V(G)-S-V(C_1))-\sum_{i=2}^q E(C_i)\\
 &\geq \ex(2m+1,K_{r+1})-sn-\ex(\sqrt{\rho}n,K_{r+1})\\
 &>\ex(2m+1,K_{r+1})-\rho^{2/}n^2.
\end{align*}
Since $G[S\cup V(C_1)]$ is $K_{r+1}$-free, by Theorem \ref{Kr-free-close},
$G[S\cup V(C_1)]$ is $\rho^{1/3}$-close to $T(|S\cup V(C_i)|,K_{r+1})$. It follows that
$G$ is $\rho^{1/4}$-close to $T(2m+1,K_{r+1})\cup (n-2m-1)K_1$.

\medskip
\textbf{Case 2.~}$s\geq m-8\rho n$.
\medskip

Let $G':=G-(\cup_{i\in [q]}E(C_i))$.  Recall $\sum_{i=1}^q (|C_i|-1) \leq 2(m-s)$ by \eqref{para} and Lemma \ref{combturan}. On the one hand, we have
\begin{align*}
e(G')&>e(G)-\sum_{i=1}^q \ex(|V(C_i)|,K_{r+1})\\
&\geq e(G)-\ex(\sum_{i=1}^q (|C_i|-1)+1,K_{r+1})\\
&\geq e(G)-\ex(2(m-s)+1,K_{r+1})\\
&\geq g(n,m,r)- 120 \rho^2 n^2  \\
&\geq (1-1/r)m^2/2+(n-m)m- 129 \rho^2n^2.
\end{align*}
On the other hand, we have
\begin{align*}
e(G')= e_G(\cup_{i\in [q]}V(C_i),S)+e(G[S]).
\end{align*}
 Note that
\begin{align}\label{G[S-D-up]}
e_G(\cup_{i\in [q]}V(C_i),S)\leq (n-s)(s-b)\leq (n-m)(m-b)
\end{align}
and
\begin{align*}
e(G[S])&\leq \ex(s-b,K_r)+b(s-b) \\
&\leq \ex(m-b,K_r)+b(m-b) \\
&\leq \ex(m,K_r)+b(m-b).
\end{align*}
That is
\begin{equation} \label{e(G')-up}
e(G[S])\leq (1-1/r)m^2/2+b(m-b).
\end{equation}

From (\ref{e(G')-up}), (\ref{G[S-D-up]}) and the lower bound for $e(G')$,  we infer that $b\leq \rho^{2/3}n$.
Furthermore, we get
\begin{align}\label{bi-close}
e_G((\cup_{i\in [q]}V(C_i),S)\geq m(n-m)-\rho^{1/2}n^2.
\end{align}
We choose $x\in S$ with $d_{G[S]}(x)=\Delta(G[S])$ and  $W\in {S\choose b}$ such that $G[S]-W$ is $K_r$-free.
Let
\[
A=\left\{
  \begin{array}{ll}
    S-W, & \hbox{if $\Delta(G[S])\leq |S|-|W|$;} \\
    N_{G[S]}(x), & \hbox{otherwise.}
  \end{array}
\right.
\]
From (\ref{e(G')-up}) and (\ref{G[S-D-up]}), we have
\begin{align*}
e(G[A])&\geq e(G[S])-bm\\
&\geq (1-1/r)m^2/2-\rho^{1/2} m^2.
\end{align*}
Thus by Theorem \ref{Kr-free-close}, $G[A]$ is $\rho^{1/2}$-close to $T(|A|,K_r)$.
By (\ref{bi-close}), $G$ is  $\rho^{1/3}$-close to $G(n,m,r)$. This completes the proof. \qed

\subsection{Perfect fractional matching}

We now prove the main result of this section. A \emph{fractional matching} in a $k$-graph $H$ is a function $w:E\rightarrow[0,1]$ such that for any $v\in V(H)$, $\sum_{e\in E:v\in e}w(e)\leq 1$. A fractional matching $w$ is \emph{perfect} if $\sum_{e\in E}w(e)=|V(H)|/k$.



Let $H$ be a balanced $(1,2)$-partite 3-graph  with vertex set $Q,P$, where $|Q|=n/2$ and $|P|=n$. For two  sets  $S_1\subset Q$, $S_2\subset P$, let $e_H(S_1,S_2)$ be the number of edges with one vertex from $S_1$ and two vertices from $S_2$. For three pairwise disjoint sets $T_1\subset Q$, $T_2, T_3\subset P$, let $e_H(T_1,T_2,T_3)$ denote the number
of edges containing exactly one vertex
from each of $T_1,T_2,T_3$.

\begin{lemma}\label{frac-PM}
Let  $\beta_1,\beta_2$ be two constants and let $n,m,t$ be three integers such that $0<1/n\ll \beta_1\ll \beta_2\ll c\ll 1$ and $t=\frac{n-2m-2}{2}\geq cn$.  Let $H$ be a balanced $(1,2)$-partite 3-graph with vertex set $Q,[n]$, where $Q=\{x_1,\ldots,x_{n/2}\}$.  Suppose $\{Q_1,Q_2\}$ is a partition of $Q$ such that $Q_1:=\{v_1,\ldots,v_t\}$ and $Q_2:=Q\setminus Q_1$. If the following five conditions hold, then $H$ has a perfect fractional  matching.
\begin{itemize}
  \item [(i)] For every $x\in Q_1$, $|N_H(x)|\geq {n\choose 2}-\beta_1 n^2$;
  \item [(ii)]for every $x\in Q_2$, $|N_H(x)|\geq f(n,m,r)-\beta_1 n^2$;
 \item[(iii)] for every vertex $y\in [n]$, $d_{H-Q_2}(y)\geq t{n\choose 2}-\beta_1 n^2$;
  \item [(iv)] for any $U\in {[n]\choose n-m}$, $e_{H-Q_1}(Q_2,U)\geq \beta_2 n^3$;
   \item [(v)] for any $U\in {[n]\choose 2m+1}$, $e_{H-Q_1}(Q_2,U,[n]-U)\geq \beta_2 n^3$.
\end{itemize}

\end{lemma}

\pf Let $\beta_1'$ be a constant such that $\beta_1\ll \beta_1'\ll \beta_2$. Let $\omega:V(H)\rightarrow \mathbb{R}^+\cup \{0\}$ be a minimum
fractional vertex cover of $H$.
Suppose that $\omega(v_{t+1})\geq \cdots\geq\omega(v_{n/2})$ and $\omega(1)\geq \cdots \geq \omega(n)$.
Let $H'$ be the $(1,2)$-partite $3$-graph with vertex set $V(H)$ and edge set $E(H')=E'$, where
\[
E'=\left\{e\in {V(H)\choose 3}\ : \ |e\cap Q|=1 \ \mbox{and }\sum_{x\in e}\omega(x)\geq 1\right\}.
\]
By the definition of $H'$, we have $E(H)\subseteq E(H')$ and $\omega$ is a fractional vertex cover of $H'$. Note that every fractional vertex cover of $H'$ is also a vertex cover of $H$. Thus $\omega$ is also a minimum fractional vertex cover of $H'$.
Let $\nu_f(H)$ and $\nu_f(H')$ denote the maximum fractional matching numbers of $H$ and $H'$, respectively;
then by the Strong Duality Theorem of linear programming, $\nu_f(H)=\omega(H)$ and  $\nu_f(H')=\omega(H')$.
Thus $\nu_f(H)=\nu_f(H')$ and hence, it suffices to show that $H'$ has a perfect matching.

Next, we observe that the edges of $H'$ form a stable family with respect to the above ordering of vertices in $P$ and $Q$:
for any  $e_1=\{v_{i_1},i_2,i_3\}$ and $e_2=\{v_{j_1},j_2,j_3\}$ with
$i_\ell\leq j_\ell$ for $1\leq \ell\leq 3$, $e_2\in E(H')$ implies  $e_1\in E(H')$.
To see this, note that, since $i_\ell\leq j_\ell$ for $1\leq \ell\leq 3$,
we have $\omega(v_{i_1})\geq \omega(v_{j_1})$ and $\omega(i_\ell)\geq \omega(j_\ell)$ for $2\leq \ell\leq 3$.
If $e_2\in E(H')$ then $ \sum_{x\in e_2}\omega(x)\geq 1$. We get $\sum_{x\in e_1}\omega(x)\geq 1$ and thus  $e_1\in E(H')$.

Let $G$ be a graph with vertex set $[n]$ and edge set $N_{H'}(x_{n/2})$. Then $G$ is a stable 2-graph.

\medskip
\textbf{Claim 1.~} $\nu(G)\leq m$ .
\medskip

Suppose, to the contrary that $\nu(G)\geq m+1$. Let $M_1:=\{e_1,\ldots,e_{m+1}\}$ be a matching of $G$.  Let $M_1':=\{e_i\cup \{v_{i+t}\}\ |\ i\in [m+1]\}\}$. Note that for any $x\in Q_1$, we have $d_{H-V(M_1')}(x)\geq \binom{n-2m-2}{2})-\beta_1 n^2$ by assumption (i). In addition, assumption (iii) implies that for $y\in [n]-V(M_1')$, we have $d_{H-V(M_1')}\geq t\binom{n-2m-3}{2}-\beta_1 n^2$.
Let $\mathcal{A}(n-2m-2)$ denote the $(1,2)$-paritite graph with biaprtition $(Q_1,[n]-V(M_1))$ such that for any $x\in Q_1$, $N_{\mathcal{A}(n-2m-2)}(x)={[n]-V(M_1)\choose 2}$.
Then one can see that every vertex $v\in V(H-V(M_1'))$ is $\beta_2$-good with respect to $\mathcal{A}(n-2m-2)$. So by Lemma \ref{good-1}, $H-V(M_1)$ has a matching $M_2$ of size $t$.
 Then $M_1'\cup M_2$ is a perfect matching of $H$. This completes the proof of Claim 1.

\medskip
\textbf{Claim 2.~} $G$ is  $\beta_1'$-close to $S(n,m)$ or $D(n,m)$.
\medskip

Otherwise, since $d_{H'}(v_{n/2})\geq f(n,m,r)-\beta_1n^2$, then by Lemma \ref{stable-lem}, $G$ has a matching $M_1$ of size $m+1$, contradicting to Claim 1.
 This completes the proof of Claim 2.

Next we discuss two cases.

\medskip
\textbf{Case 1.~} $G$ is  $\beta_1'$-close to $S(n,m)$.
\medskip

Recall that $G$ is stable. We claim that for $1\leq i\leq m-\sqrt{\beta_1'}n$, $|N_G(i)\cap ([n]-[m])|> n-m-\sqrt{\beta_1'}n$. Otherwise, suppose that $|N_G(\sqrt{\beta_1'}n)\cap ([n]-[m])|\leq n-m-\sqrt{\beta_1'}n$. Then  we have
\[
|E(S(n,m))\backslash E(G)|> \beta_1' n^2,
\]
a contradiction.

By (iv), we have $e_{H-Q_1}(Q_2,[n]-[m])\geq \beta_2n^3$. Let $b:=\sqrt{\beta_1'}n+1$. So we may greedily find a matching $M_1$ of size $b$ in $H-Q_1-[m]$.

Let $G':=G-V(M_1)$. Next we greedily find a matching of size $m-b+1$ in $G'$. Since $e_G(\{b\},[n]-[m])\geq n-m-b$, then we find $v_b\in [n]-[m]-V(M_1)$ such that $e_1:=bv_b\in E(G)$. Suppose we have constructed $\ell$ pairwise disjoint edges $e_1,\ldots,e_\ell$ such that $\ell<m-b+1$, for $1\leq i\leq m-b+1$, $m-b+2-i\in e_i$ and $e_i-\{b+1-i\}\subseteq [n]-[m]-V(M_1)$.
Since $e_G(\{m-b+2-\ell\},[n]-([m]\cup V(M_1))\geq n-m-\ell-1-b\geq cn$, we may find $v_{\ell+1}\in [n]-([m]\cup V(M_1)\cup (\cup_{i=1}^\ell e_i))$ such that $\ell'v_{\ell+1}\in E(G)$, where $\ell':=m-b+2-\ell$. Continuing this process, we may find a matching $M_2'$ of size $m-b+1$ in $G'$. Let $M_2':=\{f_1,\ldots,f_{m-b+1}\}$ and let $Q_2-V(M_1)=\{x_{i_1},\ldots,x_{i_{m-b+1}}\}$. Define $M_2:=\{f_j\cup \{x_{i_j}\ |\ 1\leq j\leq m-b+1\}\}$.
Then $M_1\cup M_2$ is a matching of size $m+1$ in $H-Q_1$.

Let $\mathcal{A}(t)$ be a complete $(1,2)$-partite $3$-graph with partition $(Q_1, [n]-V(M_1\cup M_2)$
Let $H_1:=H-V(M_1\cup M_2)$. By (i) and (iii), we know  that every vertex of $H_1$ is $\sqrt{\beta_2}$-good with respect to $\mathcal{A}(t)$. Thus by Lemma \ref{good-1}, $H_1$ has a matching $M_3$ of size $t$ and $M_1\cup M_2\cup M_3$ is a perfect matching  of $H$.

\medskip
\textbf{Case 2.~} $G$ is  $\beta_1'$-close to $D(n,m)$.
\medskip

 Since $G$ is $\beta_1'$-close to $D(n,m)$, all but at most
 $2\sqrt{\beta_1'}n$ vertices of $G$ are $\sqrt{\beta_1'}$-good with respect to $D(n,m)$. Without loss of generality, suppose that $U=[2m+1]$ and for $x\in [n]-[2m+1]$, $d_{D(n,m)}(x)=0$.
Let $U^{bad}$ denote the set of $\sqrt{\alpha_0}$-bad vertices in  $U$. Let $h:=|U^{bad}|$ and  $U_1=U-U^{bad}$. Without loss of generality, suppose  that $|U_1|$ is even.
Since $G$ is stable, we may assume that $U^{bad}=\{2m-h+2,\ldots,2m+1\}$.
By (v), $e_{H-Q_1}(Q_2,[n]-U_1,[U_1])\geq \beta_2 n^3$.
So we may greedily find a matching $M_1$ of size $h+1$. Let $G':=G[[2m+1]]-V(M_1)$. We know that for every $y\in V(G')$,
\[
d_{G'}(y)\geq (1-1/r)m-6\sqrt{\beta_1'}n\geq |V(G')|/2.
\]
Thus $G'$ has a perfect matching $M_2'$ of size $|U_1|/2$.
Let $M_2':=\{f_1,\ldots,f_{|U_1|/2}\}$ and $T:=\{v_{i_1},\ldots,v_{i_h}\}=U_1-V(M_2')$.  Define $M_2:=\{f_j\cup \{v_{i_j}\}\ | j\in [|U_1|/2]\}$
Thus $M_1\cup M_2$ is a matching of size $m+1$ in $H-Q_1$. By the same argument as in Case 1, $H-Q_1-V(M_1\cup M_2)$ has a perfect matching $M_3$ of size $t$. Then $M_1\cup M_2\cup M_3$ is a perfect matching of $H$.
This completes the proof. \qed

\subsection{Almost perfect matching}

In this section, we aim to find a ``near regular" spanning subgraph of $H$.
The discussion here follows that in \cite{AFH12}.
We need to find a   random subgraph of a balanced $(1,2)$-partite 3-graph
 and use it to find a subgraph on which a ``R\"odl nibble" result can be applied.

We also need the following result attributed to Pippenger  and Spencer \cite{NJ} to find an almost perfect matching in $H$. A {\it cover} in a hypergraph $H$
is a set of edges whose union is $V(H)$.

\begin{lemma}[Pippenger and Spencer, 1989] \label{nibble}
	For every integer $k\ge 2$ and reals $r\ge 1$ and $a>0$, there are $\gamma=\gamma(k,r,a)>0$ and $d_0=d_0(k,r,a)$ such that for every $n$ and $D\ge d_0$ the following holds: Every $k$-uniform hypergraph $H=(V,E)$ on a set $V$ of $n$ vertices in which all vertices have positive degrees and which satisfies the following conditions:
	\begin{itemize}
	 	\setlength{\itemsep}{0pt}
		\setlength{\parsep}{0pt}
		\setlength{\parskip}{0pt}
		\item[$($1$)$] For all vertices $x\in V$ but at most $\gamma n$ of them, $d(x)=(1\pm \gamma)D$;
		\item[$($2$)$] For all $x\in V$, $d(x)<rD$;
		\item[$($3$)$] For any two distinct $x,y\in V$, $d(x,y)<\gamma D$;
	\end{itemize}
	contains a cover of at most $(1+a)(n/k)$ edges.
\end{lemma}

Note that $H$ contains a cover of at most $(1+a)(n/k)$ edges implies
that $H$ contains a matching of size at least $(1-(k-1)a)(n/k)$ (see, for example, \cite{NJ}).

\begin{theorem}[F\"{u}redi, \cite{Fured}]\label{Fu81}
For any $k$-uniform hypergraph $H$ with $n$ vertices,  there exists a sub-hypergraph $H'$ of $H$ such that $e(H')\leq n$ and $\nu_f(H)=\nu_f(H')$.
\end{theorem}

\begin{lemma}\label{Near-PM}
Let $n,r,m,t$ be four positive integers  and let $\beta,c$ be two constants such that $n$ is even, $n\geq 2m+cn$ and $1/n\ll \gamma\ll \beta\ll c\ll 1/r$. Let $Q_1,Q_2$ be two vertex disjoint subsets such that $Q_1=t$ and $Q_2:=n/2-t$. Let $Q=Q_1\cup Q_2$. Let $H$ be a $(1,2)$-partite 3-graph with vertex partition $(Q, [n])$. If the following four conditions, then $H$ has a matching covering all but at most $\gamma n$ vertices.
\begin{itemize}
  \item [(i)] For every $x\in Q_1$, $|N_H(x)|={[n]\choose 2}$;
  \item [(ii)]for every $x\in Q_2$, $|N_H(x)|\geq f(n,m,r)-\beta_1 n^2$;
  \item [(iii)] for any $U\in {[n]\choose n-m}$, $e_{H-Q_1}(Q_2,U)\geq \beta_2 n^3$;
   \item [(iv)] for any $U\in {[n]\choose 2m+1}$, $e_{H-Q_1}(Q_2,U,[n]-U)\geq \beta_2 n^3$.
\end{itemize}

\end{lemma}

\pf Firstly, we establish the following claim.

\medskip
\textbf{Claim 1.}~$H$ has $n/\ln n$ edge disjoint perfect fractional  matchings $g_1,\ldots,g_{n/\ln n}$ such that for any $x,y\in V(H)$,
$\sum_{i=1}^{n/\ln n}\sum_{\{x,y\}\subseteq e\in E(H)}g_i(e)\leq 3$.

By Lemma \ref{frac-PM}, $H$ has a perfect fractional  matching. Then by Lemma \ref{Fu81}, $H$ has a perfect fractional  matching $g_1$ such that $|E_1|\leq 3n/2$,
where
\[
E_1:=\{e\in E(H_1)\ |\ g_1(e)>0\}.
\]
Let $H_0:=H$ and $H_1:=H_0-E_1$.
Suppose we have constructed $s$ graphs $H_1,\ldots, H_s$ and  $s$ edge-disjoint perfect fractional  matchings $g_1,\ldots,g_s$ such that for $1\leq i\leq s$, $g_i$ is a perfect fractional  matching of $H_{i+1}$ for any $x,y\in V(H)$, $\sum_{i=1}^s\sum_{\{x,y\}\subseteq e}g_i(e)\leq 3$.
Let $$\mathcal{A}_s:=\{\{x,y\}\in {V(H)\choose 2}\ |\ \sum_{i=1}^s\sum_{\{x,y\}\subseteq e}g_i(e)\geq 2\}.$$
Let
\[
E_{s_1}:=\{e\in E(H)\ |\ \exists e'\in \mathcal{A}_S\mbox{ such that $e'\subseteq e$}\},
\]
and
\[
E_{s_2}:=\{e\in E(H_{s-1})\ |\ g_s(e)>0\}.
\]
Define $E_s:=E_{s_1}\cup E_{s_2}$ and $H_{s}:=H_{s-1}-E_s$. Let $G_s$ denote the graph with vertex set $V(H)$ and edge set $\mathcal{A}_s$. Observe that for every vertex $x\in V(H)$, $\sum_{i=1}^s\sum_{x\in e}g_i(e)=s$, which implies  $d_{G}(x)\leq s/2$ for all $x\in V(H)$. For any $x\in Q$, we then obtain
\begin{align*}
d_{H_s}(x)&\geq d_{H}(x)-d_{G_s}(x)n-\sum_{i=1}^s|E_s|\\
&\geq d_H(x)-2sn.
\end{align*}
Consequently, $H_{s}$ still satisfies the four conditions in Lemma \ref{frac-PM}, so $H_s$ admits a perfect fractional  matching $g_s$. Iterating this procedure, we obtain the desired collection of $n/\ln n$ edge‑disjoint perfect fractional matchings $g_1,\dots,g_{n/\ln n}$. This completes the proof of Claim 1.

By Claim 1, let $g_1,\dots,g_{n/\ln n}$ be edge‑disjoint perfect fractional matchings satisfying
\[
\sum_{i=1}^{n/\ln n}\sum_{\{x,y\}\subseteq e\in E(H)}g_i(e)\leq 3
\]
for all $x,y\in V(H)$.
Define $h(e):=\sum_{i=1}^{n/\ln n}g_i(e)$.
Since the supports $E(g_1),\dots,E(g_{n/\ln n})$ are pairwise disjoint, we have $0\le h(e)\le 1$ for every edge $e\in E(H)$.

We construct a random spanning subgraph $F\subseteq H$ by including each edge $e$ independently with probability $h(e)$.
For each vertex $v\in V(H)$,
\[
\mathbb{E}[d_F(v)]=\sum_{e\ni v}h(e)=\frac{n}{\ln n},
\]
as each $g_i$ is a perfect fractional matching.
For every $2$-element set $D\subseteq V(H)$,
\[
\mathbb{E}[d_F(D)]=\sum_{e\supseteq D}h(e)\le 3.
\]
Applying the Chernoff bound, with probability $1-o(1)$ there exists a spanning subgraph $F'$ such that
\[
\bigl(1-1/\log n\bigr)\frac{n}{\ln n} \le d_{F'}(v)\le \bigl(1+1/\log n\bigr)\frac{n}{\ln n}
\]
for all $v\in V(H)$, and furthermore $\Delta_2(F')\le \sqrt{n}$.
Invoking Lemma \ref{nibble}, we conclude that $H$ contains the desired matching.
This finishes the proof.\qed

\subsection{Proof of Theorem \ref{main}}

\textbf{Proof of Theorem \ref{main}} Given a system of graphs $\G=\{G_1,\ldots,G_{m+1}\}$ satisfying the assumptions in Theorem \ref{main}, we need only to consider the case where $2m+cn \leq  n \leq 10m$ for some constant $c>0$.
	By Lemmas~\ref{Phk} and \ref{Divesity-close}, we may assume $H(\cal G)$ is not
        $\varepsilon$-close to  $H_S(n,m)$ or $H_D(n,m)$, where $\varepsilon\ll 1$. Now we consider the $(1,2)$-partite 3-graph $H^{\ast}(\cal G)$. To simplify the notation, we next write $H^{\ast}$ for $H^{\ast}(\cal G)$.
		 Lemma~\ref{absorbing} implies that $H^{\ast}$ contains  a matching $M'$ such that, for some $0<\rho'\ll \rho \ll \varepsilon$,
      $|M'|\leq \rho n/4$ and,
     for any  $S\subseteq V(H^\ast)$ with $|S|\leq \rho' n$
     and $2|S\cap Q|=|S\cap P|$, $H^\ast[S\cup V(M)]$ has a perfect matching.

	  Let $H_1^{\ast}= H^\ast-V(M')$. Then $d_{H^{\ast}_1}(v)>\binom{n'}{2}-\rho (n')^2,\ d_{H_1^{\ast}}(\{v,u\})>\beta_2n'-\rho n'$ for all $v\in Q_1\cap V(H_1^\ast)$ and $u\in P\cap V(H_1^\ast)$, and $d_{H_1^\ast}(v)>f(n',m,r)-\rho (n')^2$ for all $v\in Q_2\cap V(H_1^\ast)$. Furthermore,
      $H_1^\ast-Q_1$ is not $(2\varepsilon)$-close to $H_D(n',m')$ or $H_S(n',m')$ by Lemmas~\ref{Phk} and \ref{Divesity-close}, where $n'=n-2|M|'$ and $m=m-|M'|$.

	By Lemma~\ref{Near-PM}, $H_1^\ast$ contains a matching $M_1$ covering
        all but at most $\gamma n$ vertices, where we
    choose $\gamma$ so that $0<\gamma<\rho'$.
	Now	$H[(V(H_1^\ast)\setminus V(M_1))\cup V(M')]$ has a perfect matching $M_2$. Clearly, $M_1\cup M_2$ forms a perfect matching in $H^\ast$.\qed

\end{document}